\def\comment#1{{}}

\documentclass[11pt]{amsart}
\usepackage{amssymb, amsmath, amscd, dcpic, pictexwd}
\usepackage{hyperref} 

\numberwithin{equation}{section}
\newtheorem{thm}{Theorem}[section]

\newtheorem{cor}[thm]{Corollary}
\newtheorem{lem}[thm]{Lemma}
\newtheorem{prop}[thm]{Proposition}
\newtheorem{defn}[thm]{Definition}
\newtheorem{rem}[thm]{Remark}
\newtheorem{exm}[thm]{Example}

\newcommand{\cA}{\mathcal A}

\newcommand{\cD}{\mathcal D}

\newcommand{\cS}{\mathcal S}
\newcommand{\cO}{\mathcal{O}}
\newcommand{\cL}{\mathcal{L}}

\newcommand{\cE}{\mathcal E}

\newcommand{\cM}{\mathcal M}
\newcommand{\cH}{\mathcal H}

\newcommand{\cY}{\mathcal Y}
\newcommand{\R}{{\bf R}}

\newcommand{\fb}{\mathfrak{b}}

\newcommand{\fg}{\mathfrak{g}}
\newcommand{\fh}{\mathfrak{h}}
\newcommand{\gh}{\fh}

\newcommand{\fl}{\mathfrak{l}}

\newcommand{\fp}{\mathfrak{p}}

\newcommand{\fs}{\mathfrak{s}}
\newcommand{\ft}{\mathfrak{t}}

\newcommand{\HH}{\mathbb{H}}

\newcommand{\CC}{\mathbb{C}}
\newcommand{\RR}{\mathbb{R}}
\newcommand{\C}{\CC}
\newcommand{\PP}{\mathbb{P}}
\newcommand{\ZZ}{\mathbb Z}
\newcommand{\Z}{\ZZ}
\newcommand{\into}{\hookrightarrow}
\newcommand{\onto}{\twoheadrightarrow}

\newcommand{\ra}{\rightarrow}
\newcommand{\half}{{1\over 2}}
\newcommand{\bra}{{\langle}}
\newcommand{\ket}{{\rangle}}

\newcommand{\br}{\buildrel}
\newcommand{\bs}{\bigskip}
\newcommand{\cut}{\vskip-.3in}
\newcommand{\blank}{\hskip.3in}

\newcommand{\Sym}{{\mbox {Sym~}}}
\newcommand{\End}{{\mbox {End~}}}
\newcommand{\Hom}{{\mbox{Hom}}}
\newcommand{\pr}{{\mbox{pr}}}
\newcommand{\Aut}{{\mbox{Aut~}}}
\newcommand{\Det}{{\mbox{det}}}
\newcommand{\Vect}{{\mbox{Vect}}}
\newcommand{\Pic}{{\mbox{Pic}}}
\newcommand{\Res}{{\mbox{Res~}}}
\newcommand{\id}{{\mbox{id}}}
\newcommand{\Der}{{\mbox{Der~}}}
\newcommand{\ts}{\textsection}

\author{Bong H. Lian and Shing-Tung Yau}

\title{Period Integrals of CY and General Type Complete Intersections ${}^1$\footnote{${}^1${A\lowercase{ppeared in }I\lowercase{nvent.} M\lowercase{ath.} V\lowercase{ol.} 191, 1 (2013) 35-89. }}}

\begin{document}
\maketitle
\begin{abstract}
We develop a global Poincar\'e residue formula to study period integrals of families of complex manifolds. For any compact complex manifold $X$ equipped with a linear system $V^*$ of generically smooth CY hypersurfaces, the formula expresses period integrals in terms of a canonical global meromorphic top form on $X$. Two important ingredients of this construction are the notion of a CY principal bundle, and a classification of such rank one bundles.
We also generalize the construction to CY and general type complete intersections. When $X$ is an algebraic manifold having a sufficiently large automorphism group $G$ and $V^*$ is a linear representation of $G$, we construct a holonomic D-module that governs the period integrals. The construction is based in part on the theory of tautological systems we have developed in the paper \cite{LSY1}, joint with R. Song. The approach allows us to explicitly describe a Picard-Fuchs type system for complete intersection varieties of general types, as well as CY, in any Fano variety, and in a homogeneous space in particular. In addition, the approach provides a new perspective of old examples such as CY complete intersections in a toric variety or partial flag variety.

\end{abstract}

\tableofcontents
\baselineskip=16pt plus 1pt minus 1pt
\parskip=\baselineskip

\pagenumbering{arabic}
\addtocounter{page}{0}
\markboth{\SMALL Bong H. Lian and Shing-Tung Yau}
{\SMALL Period Integrals of CY and General Type Complete Intersections}

\section{Introduction}

Let $\pi:\cY\ra B$ be a family of $n$-dimensional CY manifolds. Recall that $R^n\pi_*\C$ is a flat vector bundle on $B$ whose fiber is the middle cohomology $H^n(Y_b,\C)$, where $Y_b=\pi^{-1}(b)$. Thus, assuming that $B$ is simply connected, fixing a base point $b_0\in B$ gives us a canonical trivialization of this bundle with fixed fiber $H^n(Y_{b_0},\C)$. The bundle contains a subbundle $\HH^{top}$ whose fiber at $b\in B$ is $H^0(K_{Y_b})$. Consider a nonzero local section $\omega$ of $\HH^{top}$ over a small neighborhood $U$ of $b_0$. The  local Torelli theorem implies that the line $\C\omega_b=H^0(K_{Y_b})\subset H^n(Y_{b_0},\C)$ determines the isomorphism class of $Y_b$. One way to study the variation of this line is to consider the period integrals
\begin{equation}\label{local-periods}
\int_\gamma\omega_b,\blank \gamma\in H_n(Y_{b_0},\Z)
\end{equation}
as functions defined on $U$, attempt to construct differential equations for them, and ultimately compute them. This presents two problems. First, the normalization of the functions \ref{local-periods} depends on the choice of local section $\omega$, rendering any system of differential equations for \ref{local-periods} dependent on such an arbitrary choice. Second, any description of such a system of differential equations is inherently local, depending on an arbitrary choice of $U$. 
As will be seen shortly, we can circumvent both of these problems by introducing a canonical global normalization for the period integrals, in the case when $\cY$ is a family of complete intersections in a fixed ambient manifold $X$. 

We begin with a summary of the main results and ideas.
Let $X$ be a compact connected $d$-dimensional complex manifold, and $L_1,..,L_s$ be line bundles on $X$ such that
$$
V_i:=H^0(L_i)^*\neq 0
$$
and that the general section $\sigma_i\in V_i^*$ defines a nonsingular hypersurface $Y_{\sigma_i}=\{\sigma_i=0\}$ in $X$. 
We would like to consider the family of all complete intersections 
$$
Y_\sigma:=Y_{\sigma_1}\cap\cdots\cap Y_{\sigma_s}
$$
which are smooth of codimension $s<d$ in $X$. Put
$$
V:=V_1\times\cdots\times V_s,\blank B=V^*-D
$$
where $D$ consists of $\sigma=(\sigma_1,..,\sigma_s)\in V^*$  such that $Y_\sigma$ is not smooth of codimension $s$. Then $B$ parameterizes a smooth family $\cY$ of smooth complete intersections of codimension $s$ in $X$. Given $\sigma\in B$, the adjunction formula gives a canonical isomorphism
$$
(L+K_X)|Y_\sigma\cong K_{Y_\sigma}
$$ 
where $L:=L_1+\cdots+L_s$. Thus when $L+K_X$ is trivial, $\cY$ is a family of CY manifolds, and when $L+K_X$ is ample, $\cY$ is a family of general type manifolds.

We shall assume that the Hodge number $\dim H^0(K_{Y_\sigma})$ is a locally constant function on $B$ (this is true when the fibers of $\cY$ are K\"ahler), and consider bundle $\HH^{top}$ whose fiber at $\sigma\in B$ is $H^0(K_{Y_\sigma})$. For $\sigma\in B$, let 
$$\R_\sigma:H^0(L+K_X)\ra H^0(K_{Y_\sigma})$$ 
be the restriction map. Clearly, for each $\tau\in H^0(L+K_X)$, the map
$$
B\ra \HH^{top},~~~\sigma\mapsto \R_\sigma(\tau)
$$
defines a global section of $\HH^{top}$. So, we get a linear map
$$
\R:H^0(L+K_X)\ra H^0(B,\HH^{top}),~~~\tau\mapsto(\sigma\mapsto\R_\sigma(\tau)).
$$
In particular when $L+K_X=\cO_X$, 
then $\R(1)$ gives a canonical global trivialization of $\HH^{top}$. 

\begin{defn} (Period sheaves)\label{period-integrals}
We call the map $\R$ the global Poincar\'e residue map of the family $\cY$. For each $\tau\in H^0(L+K_X)$, we define
the period sheaf ${\bf\Pi(\tau)}$ of the family $\cY$ to be the locally constant sheaf on $B$ generated by the local sections
$$
\int_\gamma\R(\tau),\blank\gamma\in H_{d-s}(Y_\bullet,\ZZ)
$$
where $Y_\bullet$ is some fixed fiber of $\cY$. A local section of this sheaf is called a period integral.
\end{defn}

The main goal of this paper is to give a new description of the map $\R$, and use it to construct explicitly a system of Picard-Fuchs type differential equations that govern the period integrals. The resulting system will turn out to be a certain generalization of a {\it tautological system}. The latter notion was introduced in \cite{LSY1}, where it was applied to the special case when $X$ is a partial flag variety and $L+K_X$ is trivial. In this special case, period integrals were defined in an ad hoc fashion, by choosing a particular meromorphic form to write down a Poincar\'e residue. The resulting system of differential equations was shown to be holonomic, and was amenable to fairly simple description. However, the construction was clearly specific to that example and to the case of CY complete intersections, and it gives $\R$ only up to an overall normalization by an undetermined holomorphic function on the base $B$. In this paper, we give a new construction of $\R$ that removes the normalization ambiguity entirely, and at the same time, generalizes to an arbitrary manifold $X$ and to a general family $\cY$, not necessarily CY. 

The first problem is to describe the map induced on global sections by the Poincar\'e residue map  $PR:\Omega_X^d(Y)\ra\Omega^{d-s}_Y$, in a way that is compatible with the variation of $Y$. The map $PR$ is a sheaf homomorphism, whose only known description (except in isolated cases), unfortunately, relies on local coordinates. The variation of $Y$ in a family also makes matter worse. Since period integrals are inherently global objects on $Y$ (and global also on the parameter space if we lift to its universal cover), it is difficult to work with the local description of $PR$ for the purpose of constructing and describing their differential equations, let alone computing the period integrals. A key insight in our approach is that the map $\R$ can be realized as a family version of the $PR$, but in terms of an explicit globally defined meromorphic form on the ambient space $X$. This realization is carried out by lifting the construction of period integrals to a certain principal bundle over $X$. 

One of the main benefits of this approach is that sections of line bundles can then be represented as globally defined functions, affording a description of the Poincar\'e residue map in purely global terms (Theorem \ref{globalPR}.) Quite remarkably, the result  is in fact independent of the choice of principal bundle (Theorem \ref{unique-globalPR}.) So, every choice gives a potentially new way to express the same period integrals. More importantly, one can show that such a principal bundle exists for any complex manifold $X$, and that there is always a canonical choice (Theorem \ref{Calabi-form}.) This part uses an old idea of E. Calabi. In specific examples, the global realization of $\R$ often allows us to compute power series expressions for period integrals. As it is well-known, the latter are central to mirror symmetry.

Finally, our realization of $\R$  allows us to show that the period integrals are governed by a tautological system (Theorem \ref{general-tautological}), or a certain enhanced version of it (Theorem \ref{enhanced-tautological}.) In the projective case, we apply our results to toric manifolds  and homogeneous spaces, as special examples, and show that the period sheaves are governed by {\it holonomic} tautological systems. Even in the general setting, these systems are still amenable to fairly explicit and simple descriptions. For example, in the case of $X$ a general homogeneous manifold, we give two different descriptions (Theorems \ref{homogeneous-one}, \ref{Segre-Veronese}) by using the Borel-Weil theorem and a theorem of Kostant and Lichtenstein, and we enumerate a set of generators for the tautological system. For $X$ a toric manifold, our tautological systems turn out to be examples of GKZ hypergeometric systems and their extended versions (Example \ref{toric-revisited}.) Explicit power series solutions to these systems are also known. 

As part of the toolkit in our approach, we also develop a general form of the Euler sequence for $TX$, and a principal bundle version of the adjunction formula for $K_X$ (Theorem \ref{KX-theorem}.) More generally, a principal $H$-bundle $H- M\ra X$ enables us to describe vector bundles on $X$, and their sections in fairly simple functional terms. For example, certain bundles are described by representations of $H$; sections of those bundles can be realized as global holomorphic functions on $M$. Moreover, these descriptions are all equivariant with respect to any given automorphism of $X$, whenever the $H$-bundle $M$ itself is equivariant. In special cases such as homogeneous spaces, some of these constructions are known in representation theory.

Another important application of our general Euler sequence is to the study of {\it CY structures}. Consider a principal $H$-bundle $H-M\ra X$. A nowhere vanishing form $\omega_M\in H^0(K_M)$ is called a CY structure of $M$ if $\CC\omega_M$ is a representation of $H$. We derive an obstruction for the existence of such a structure (Theorem \ref{obstruction}), and prove that a CY structure is essentially unique, if exists (Theorem \ref{uniqueness}.) We then use these results to classify all CY structures for the rank 1 bundles. These results will be crucial to proving the global Poincar\'e residue formulas for period integrals (Theorems \ref{globalPR}, \ref{unique-globalPR}) 

\noindent {\bf Acknowledgements.}
B.H.L. is partially supported by NSF FRG grant DMS-0854965, and S.T.Y. by NSF FRG grant DMS-0804454.

\section{Adjunction formula for principal bundles}\label{sec-adjunction}

Throughout this section, let $G$ and $H$ be complex Lie groups, and $M$ a complex manifold, not necessarily compact, on which the group
$$K:=G\times H$$ 
acts holomorphically. We denote the actions by
\begin{eqnarray*}
& &G\times M\ra M,~~(g,m)\mapsto gm\cr
& &H\times M\ra M,~~(h,m)\mapsto mh^{-1}.
\end{eqnarray*}
We assume that the $H$-action on $M$ is free and proper, so that the quotient $M/H$ is a complex $G$-manifold:
$$
G\times M/H\ra M/H,~~(g,[m])\mapsto [gm].
$$
We have the projection map
$$
\pi:M\ra M/H,~~~m\mapsto[m].
$$


Let $\Vect_G(X)$ be the category of holomorphic $G$-equivariant vector bundles on $X$, and $\Vect_K(M)$ the category of holomorphic $K$-equivariant vector bundles on $M$. We have the following {\it pullback} and {\it quotient} functors:
\begin{eqnarray*}
&\pi^*&:\Vect_G(X)\ra Vect_{K}(M)\cr
&/H &:\Vect_{K}(M)\ra \Vect_G(X)
\end{eqnarray*}
where 
\begin{eqnarray*}
& &(\pi^*E)_m=(m,E_{[m]})\subset M\times E\cr
& &(F/H)_{[m]}=(\cup_{m'\in[m]}F_{m'})/H\subset F/H
\end{eqnarray*} 
for $m\in M$. Here $E\in \Vect_G(X)$ and $F\in \Vect_K(M)$.
It is straightforward to show that both functors are equivalences of categories, and that we have natural isomorphisms
\begin{eqnarray*}
(\pi^*E)/H &\cong& E\cr
\pi^*(F/H) &\cong& F.\cr
\end{eqnarray*}
Specializing to line bundles, letting $Pic_K(M)$ denote the $K$-equivariant Picard group of $M$, i.e. the group of isomorphism classes of $K$-equivariant line bundles on $M$, we get a canonical isomorphism
\begin{equation}\label{Pic-Pic}
Pic_K(M)\cong Pic_G(X).
\end{equation}

We will make use of primarily vector bundles on $M$ that are trivial, but not necessarily $K$-equivariantly. The two functors above will allow us to interpolate between representation theory of $H$ and $G$-equivariant vector bundles on $X$. 

The $G$-action on $M$ induces homomorphisms
\begin{eqnarray*}
& &G\ra\Aut TM,~~g\mapsto dg\cr
& &G\ra\Aut T^*M,~~g\mapsto \delta g^{-1}.
\end{eqnarray*}
Here $\delta g^{-1}:T^*M_m\ra T^*M_{gm}$ is the inverse transpose of the differential $dg:TM_m\ra TM_{gm}$. Likewise, the $H$-action on $M$ induces similar homomorphisms. Note that because we let $H$ act on the right on $M$, we have $dh:TM_m\ra TM_{mh^{-1}}$ and $\delta h^{-1}:T^*M_m\ra T^*M_{mh^{-1}}$.

Let $E\ra M$ be a $G$-equivariant vector bundle. Then there is an induced homomorphism
$$
\Gamma: G\ra \Aut H^0(E),~~g\mapsto\Gamma_g
$$
given by
$$
\Gamma_g\sigma=g\circ\sigma\circ g^{-1}.
$$
The same is true for a $H$-equivariant bundle. For example, we have
\begin{eqnarray*}
(\Gamma_g\sigma)_m&=&dg(\sigma_{g^{-1}m})\in TM_m\cr
(\Gamma_h\sigma)_m&=&dh(\sigma_{mh})\in TM_m
\end{eqnarray*}
for $g\in G$, $h\in H$, $m\in M$, $\sigma\in H^0(TM)$.

Let $\rho:H\ra\Aut V$ be a finite dimensional holomorphic representation of $H$. Then $H$ acts on the product $M\times V$ freely and properly, and the quotient
$$
E_\rho:=M\times_H V=(M\times V)/H
$$
is a $G$-equivariant holomorphic vector bundle of rank $\dim V$ over $M/H$. The $G$ action on $E_\rho$ is defined by
$$
G\times E_\rho\ra E_\rho,~~(g,[m,v])\mapsto[gm,v].
$$
Note that $[mh^{-1},\rho(h)v]=[m,v]$ for $h\in H$. A holomorphic function
$$
\sigma': M\ra V
$$
is said to be {\bf$\rho$-equivariant} if
$$
\sigma'(mh^{-1})=\rho(h)\sigma'(m)\blank (h\in H,~m\in M.)
$$
{\it We denote the space of such functions by $\cO(M)_\rho$.}

The following is an abstraction of a well-known fact about holomorphic sections of $G$-equivariant vector bundles on homogeneous spaces.

\begin{prop}\label{function-section}
Let $\rho:H\ra\Aut V$ be a finite dimensional holomorphic representation of $H$, and $\sigma':M\ra V$ a $\rho$-equivariant function. Then $\sigma(m):=[m,\sigma'(m)]$ ($m\in M$) defines a holomorphic section of $E_\rho$. This gives a $G$-equivariant linear isomorphism 
$$\cO(M)_\rho\ra H^0(X,E_\rho),~~\sigma'\mapsto\sigma.$$
We denote the inverse map by $\sigma\mapsto\sigma'=\sigma_M$.
\end{prop}
\begin{proof}
That $\sigma$ is a holomorphic section of $E_\rho$ is clear. To see that $\sigma'\mapsto\sigma$ above is an isomorphism, we construct the inverse map.  Given a section $\sigma\in H^0(E_\rho)$, consider the image $\sigma(M/H)\subset E_\rho$. Its pre-image $\pr^{-1}\sigma(M/H)$ under the projection $\pr:M\times V\ra E_\rho $ is the graph of a holomorphic function $\sigma':M\ra V$. It is easy to check that
$$
\sigma(m)=[m,\sigma'(m)]
$$
and that $\sigma'$ is $\rho$-equivariant. The correspondence $\sigma\mapsto\sigma'$ is the inverse map we sought.

Finally, the induced linear $G$-actions on functions $\sigma'$ and sections $\sigma$ are given by
\begin{eqnarray*}
\Gamma_g\sigma'&=&(g^{-1})^*\sigma'=\sigma'\circ g^{-1}\cr
\Gamma_g\sigma&=&g\circ\sigma\circ g^{-1}~~(g\in G.)
\end{eqnarray*}
It is then easy to check that $\Gamma_g\sigma'\mapsto\Gamma_g\sigma$ under the isomorphism above.
\end{proof}

The proposition shows in particular that to get interesting line bundle of the form $E_\rho=M\times_H\C$ with sections, $M$ must be noncompact.

\begin{cor}
If $\rho:H\ra GL_1$ is onto and $H^0(E_\rho)\neq0$, then $M$ is noncompact.
\end{cor}
\begin{proof}
By the proposition, we have a nonzero holomorphic function $\sigma'\in\cO(M)_\rho$ which is $\rho$-equivariant, hence can't be locally constant. So, $M$ is not compact.
\end{proof}

The isomorphism in the proposition can be interpreted in terms of the pullback functor as follows.

\begin{prop}
For $E\in Vect_G(X)$, 
$$
\pi^*:H^0(X,E)\ra H^0(M,\pi^*E), ~~\sigma\mapsto \id_M\times\sigma\circ\pi
$$
is an injective $G$-equivariant linear map. In particular, for a given representation $\rho:H\ra\Aut V$ of $H$ and for $E_\rho=M\times_HV$, the linear map $\pi^*$ coincides with the map 
$$H^0(X,E_\rho)\ra\cO(M)_\rho\subset H^0(M,M\times V),~~\sigma\mapsto\sigma'$$ 
(under the identification $\pi^*E_\rho\equiv M\times V$.)
\end{prop}
\begin{proof}
The injectivity of $\pi^*$ follows immediately from the surjectivity of $\pi:M\ra X$. For $\sigma\in H^0(X,E_\rho)$, since we view $\pi^*\sigma\in H^0(M,M\times V)$ as a holomorphic function $\sigma':M\ra V$, our second assertion amounts to checking the identity
$$
\sigma([m])=[m,\sigma'(m)]\blank(m\in M)
$$
which is straightforward.
\end{proof}

{\it In the rest of the paper, we shall often identify a section in $H^0(X,E)$ with the holomorphic function that represents it without explicitly saying so.}

Let $\chi\in \Hom(H,\CC^\times)$. We shall refer to $\chi$ as an $H$-character. An $H$-character can also be regarded as character of the abelian group $H/[H,H]$, and vice versa. Given an $H$-character $\chi$, let $\CC_\chi$ denote the corresponding one-dimensional representation $\CC$ of $H$. We can treat $M\times\CC_\chi$ as a $K$-equivariant line bundle. So, we have a map
$$
\Hom(H,\CC^\times)\ra Pic_K(M),~~\chi\mapsto[M\times\CC_\chi].
$$
For $\chi\in\Hom(H,\CC^\times)$, put
$$
L_\chi:=M\times_H\CC_\chi.
$$
Then $L_\chi$ is $G$-equivariant via the action
$$
g[m,1]=[gm,1].
$$
Composing with the isomorphism in eqn. \ref{Pic-Pic}, we get

\begin{prop}\label{character-linebundle}
We have a group homomorphism
$$
\Hom(H,\CC^\times)\ra \Pic_G(M/H),~~~\chi\mapsto [L_\chi].
$$
\end{prop}

In general, this map can be trivial even when $H\neq 1$. In one opposite extreme, we have

\begin{prop} \label{Popov}\cite{Popov}
Let $G$ be a complex Lie group and $H$ a closed subgroup. Then we have an isomorphism
$$\Hom(H,\CC^\times)\ra \Pic_G(G/H),~~~\chi\mapsto [L_\chi].$$
\end{prop}
\begin{proof} 
Given a $G$-equivariant line bundle $L$ over $G/H$, $H$ stabilizes the fiber $L_H$ at the coset $H$. This determines a unique $H$-character $\chi$ by which $H$ acts on $L_H$. Fix a nonzero vector $v_0\in L_H$. Then the map
$$
G\times\CC\ra L,~~(g,1)\mapsto gv_0
$$
descends to a $G$-equivariant isomorphism $L_\chi\cong L$. Suppose we have a $G$-equivariant isomorphism to the trivial bundle: $L_\chi\ra G/H\times\CC$. Then the fiber of $L_\chi$ at the coset $H$ is isomorphic, as a representation of $H$, to the trivial representation $H\times\CC$. It follows that $\chi$ is the trivial character. This shows that $\Hom(H,\CC^\times)\ra \Pic_G(M/H)$ is an isomorphism.
\end{proof}

It is important to be able to compare directly sections of a line bundle on $X$ that is realized using two different, but compatible principal bundles over $X$. Suppose we have the two principal bundles related as follows:
$$
\begin{array}{ccccc}
H_2 & - & M_2 &\br\pi_2\over\ra & X\cr
\rho\downarrow & & \pi\downarrow & & ||\cr
H_1 & - & M_1 &\br\pi_1\over\ra & X\cr
\end{array}
$$
where $\rho:H_2\ra H_1$ is a group homomorphism such that for $m_2\in M_2$, $h_2\in H_2$, 
$$
\pi(m_2h_2^{-1})=\pi(m_2)\rho(h_2)^{-1}.
$$
Now let $\chi_1$ be an $H_1$-character and put $\chi_2=\chi_1\circ\rho$. Then, clearly there is a unique isomorphism of line bundles $\tilde\pi:L_{\chi_2}\ra L_{\chi_1}$, such that $[m_2,1]\mapsto[\pi(m_2),1]$.

\begin{prop}\label{compare-bundles}
Let $H_i-M_i\ra X$, $\rho:H_2\ra H_1$, $\pi:M_2\ra M_1$, be the data as stated in the preceding paragraph.
Let $\varphi_i:\cO(M_i)_{\chi_i}\ra H^0(L_{\chi_i})$, $i=1,2$, be the canonical function-section isomorphisms given by Proposition \ref{function-section}. Then
$$
\varphi_2^{-1}\circ{\hat\pi}^{-1}\circ\varphi_1=\pi^*:\cO(M_1)_{\chi_1}\ra\cO(M_2)_{\chi_2}.
$$
where $\hat\pi:H^0(L_{\chi_2})\ra H^0(L_{\chi_1})$, $\sigma\mapsto\tilde\pi\circ\sigma$.
\end{prop}
\begin{proof}
Let $\sigma_1\in\cO(M_1)_{\chi_1}$. By Proposition \ref{function-section}, $\varphi_1\sigma_1$ is characterized by the identity
$$
(\varphi_1\sigma_1)([m_1])=[m_1,\sigma_1(m_1)]\blank(m_1\in M_1)
$$
and likewise for $\varphi_2$. We will show that
$$
\varphi_1={\hat\pi}\circ\varphi_2\circ\pi^*:\cO(M_1)_{\chi_1}\ra H^0(L_{\chi_1}).
$$
First we need to check that $\pi^*\sigma_1\in\cO(M)_{\chi_2}$. This follows from
$$
(\pi^*\sigma_1)(m_2h_2^{-1})=\sigma_1(\pi(m_2)\rho(h_2)^{-1})=\chi_1(\rho(h_2))\sigma_1(\pi(m_2)).
$$

Next, since
$$
(\varphi_2\circ\pi^*)\sigma_1=\varphi_2(\sigma_1\circ\pi)\in H^0(L_{\chi_2})
$$
it follows that
$$
\varphi_2(\sigma_1\circ\pi)([m_2])=[m_2,\sigma_1(\pi(m_2))].
$$
Applying $\tilde\pi$ to both sides, we get
$$
(\hat\pi\circ\varphi_2\circ\pi^*\sigma_1)([m_2])=\tilde\pi[m_2,\sigma_1(\pi(m_2))]=[\pi(m_2),\sigma_1(\pi(m_2))].
$$
The right side agrees with $(\varphi_1\sigma_1)([\pi(m_2)])=(\varphi_1\sigma_1)([m_2])$, as desired.
\end{proof}

Given a principal bundle $H-M{\br\pi\over\ra} X$, and a line bundle of the form $L_\chi$ where $\chi$ is an $H$-character, 
we will need to be able to compare global $\chi$-equivariant functions (sections of $L_\chi$) with the local representations of sections.
We write 
$$M|U=\pi^{-1}(U).$$

\begin{prop}\label{compare-sections}
Fix an open set $U\subset X$. Then local trivializations $\alpha:L_\chi|U=M|U\times_H\C_\chi\ra U\times\C$ of $L_\chi$ are in 1-1 correspondence with holomorphic functions $\mu_\alpha:M|U\ra\C^\times$ having the property
\begin{equation}
\mu_\alpha(mh^{-1})=\mu_\alpha(m)\chi(h)^{-1}\blank(h\in H)
\end{equation}
and under this correspondence $\alpha$ and $\mu_\alpha$ are related by
\begin{equation}
\alpha[m,1]=([m],\mu_\alpha(m))\blank(m\in M|U).
\end{equation}
Moreover, given such a local trivialization $\alpha$, if $\psi\in H^0(L_\chi)$ and $\psi_M\in\cO(M)_\chi$ is the $\chi$-equivariant function representing $\psi$, then 
$$
\psi_M(m)=\psi^\alpha([m])\mu_\alpha(m)^{-1}\blank(m\in M|U)
$$
where $\psi|U=\psi^\alpha e_\alpha$ and $e_\alpha$ is the local frame of $L_\chi$ coresponding to $\alpha$.
\end{prop}

The proof is straightforward and is left to the reader.

Next, we shall analyze the canonical bundle of $X=M/H$. As we shall see later, this bundle can be completely described in terms of certain $H$-characters (cf. Proposition \ref{Popov}.) Taking a quotient  is in some sense the dual of taking a submanifold. Thus, the idea is to find a dual version of the following adjunction formula, for principal bundles.

\begin{prop} (Adjunction formula for submanifolds)
If $Y$ is a complex $G$-submanifold of $M$, then we have an $G$-equivariant isomorphism
$$K_Y\cong K_M|Y\otimes\wedge^q N_{Y/M}$$
where $N_{Y/M}$ is the normal bundle of $Y$ in $M$, and $q=\dim M-\dim Y$.
\end{prop}

\begin{thm}\label{KX-theorem} 
 Put $X=M/H$.
\begin{itemize}
\item (Generalized Euler sequence) There is a $G$-equivariant exact sequence of vector bundles on $X$:
$$
M\times_H\fh\into(TM)/H\onto TX.
$$
\item (Adjunction formula for quotients)  There is a canonical $G$-equivariant isomorphism
$$
K_X\cong K_M/H\otimes L_{\chi_\fh}.
$$
where $\chi_\fh$ is the $H$-character of the one-dimensional representation $\wedge^q\fh$ induced by the adjoint representation of $H$.
\item (Canonical bundle) If $M$ admits a nowhere vanishing holomorphic $G$-invariant top form $\omega_M$, and an $H$-character $\chi_M$ such that $\Gamma_h\omega_M=\chi_M(h)\omega_M$, ($h\in H$) then there is a canonical $G$-equivariant isomorphism
$$
K_X\cong L_{\chi_M\chi_\fh}.
$$
\end{itemize}
\end{thm}
\begin{proof}
Consider the $G\times H$-equivariant exact sequence on $M$:
$$
DH\into TM\onto\pi^*TX
$$
where $DH\subset TM$ is the holomorphic distribution generated by the free action of $H$ on $M$. Note that $H$ acts freely and properly on all three bundles, and that we have a canonical $G$-equivariant isomorphism $(\pi^*TX)/H\cong TX$. So to prove our first assertion, it remains to show that $DH\cong M\times\fh$, $G\times H$-equivariantly.

The free $H$ action on $M$ induces the Lie algebra isomorphism 
$$
\xi:\fh\ra H^0(DH), ~~~x\mapsto\xi_x,~~~(\xi_x)_m:={d\over dt}|_{t=0}m~ exp(tx)
$$ 
which is $H$-equivariant. Here $H$ acts on $\fh$ by the adjoint representation, and on $H^0(DH)$ by the induced action on sections of the $H$-equivariant bundle $DH$. Since the $G$- and the $H$-actions on $M$ commute, $G$ acts trivially on $H^0(DH)$. So, by letting $G$ act trivially on $\fh$, $\xi$ is $G$-equivariant. The fiber $DH_m$ of $DH$ at $m$ is
$$
DH_m=\{(\xi_x)_m|x\in\fh\}.
$$
Thus
$$
\nu:M\times\fh\ra DH,~~~(m,x)\mapsto(\xi_x)_m
$$
defines a $G\times H$-equivariant isomorphism. (Note that $M\times\fh$ is not $H$-equivariantly trivial unless $H$ is abelian or $\dim H=0$.)

Now by linear algebra, the dual of the generalized Euler sequence yields a $G$-equivariant isomorphism
$$
(\wedge^tT^*M)/H\cong \wedge^d T^*X\otimes (M\times_H\wedge^q\fh^*)\cong K_X\otimes L_{\chi_\fh^{-1}}
$$
where $t=\dim M$, $d=t-q$. Tensoring both sides with $L_{\chi_\fh}$, we get
$$
K_X\cong K_M/H\otimes L_{\chi_\fh}.
$$

For the last assertion, it remains to show that
\begin{equation}\label{characters}
K_M/H\cong L_{\chi_M}.
\end{equation}
Put $\omega=\omega_M$, $\chi=\chi_M$. By assumption, $\CC\omega$ is a one dimensional representation of $H$. Define
$$
\mu:M\times\CC\omega\ra K_M,~~~(m,\omega)\mapsto\omega_m.
$$
(Note: Since $\omega$ is nowhere vanishing, this gives a canonical non-equivariant global trivialization of $K_M$.) 
We now check $G$- and $H$-equivariance separately. For $h\in H$, $m\in M$,
$$
h(m,\omega)=(mh^{-1},\chi(h)\omega){\br\mu\over\mapsto}\chi(h)\omega_{mh^{-1}}.
$$
On the other hand, we have the $H$-action $\Gamma:H\ra \Aut  H^0(K_M)$ such that
$$
(\Gamma_h\omega)_m=\delta h^{-1}\omega_{mh}.
$$
The left side is $\chi(h)\omega_m$. It follows that
$$
\delta h^{-1}\mu(m,1)=\delta h^{-1}\omega_m=(\Gamma_h\omega)_{mh^{-1}}=\chi(h)\omega_{mh^{-1}}=\mu(h(m,1)).
$$
This shows that $\mu$ is also $H$-equivariant. Similarly, it is $G$-equivariant. This gives a canonical $G$-equivariant isomorphism
$$
(M\times\CC\omega)/H\cong K_M/H.
$$
Since the left side is $L_\chi$, this yields eqn.\ref{characters}.
\end{proof}

When $X$ is a toric variety, the generalized Euler sequence were found in \cite{BatyrevCox1994}\cite{Jaczewski1994} by methods different from that used above.

\begin{defn}
We denote by $H-M\ra X$ a $G$-equivariant principal $H$-bundle over $X$. We call a nowhere vanishing form $\omega\in H^0(K_M)^G$, a CY structure of $M$, if $\CC\omega$ is a one dimensional representation of $H$, i.e. there is an $H$-character $\chi$ such that 
\begin{equation}\label{chi-omega}
\Gamma_h\omega=\chi(h)\omega\blank(h\in H.)
\end{equation}
In this case, we also call the pair $(\omega,\chi)$ a CY structure of the $H$-bundle $M$.
\end{defn}

\section{Existence and uniqueness of CY structures}\label{sec-CYstructures}

As our second application of the generalized Euler sequence, Theorem \ref{KX-theorem}, we answer the existence-uniqueness questions on CY structures. We will show that the associated bundle of $K_X$ is a CY bundle and that it is a universal one, in some sense. Throughout this section, {\it all line bundles and principal $H$-bundles $H-M\ra X$ are assumed $G$-equivariant.} 

\begin{defn}
Given a line bundle $L$ over a complex manifold $X$, let $L^\times$ denote the complement of the zero section in $L$. Equivalently, $L^\times$ is the associated principal bundle of $L$. Denote the $\C^\times$-character
$$
\chi_1:\C^\times\ra\C^\times,~~~h\mapsto h.
$$
\end{defn}

\begin{lem}\label{L-star}
For any line bundle $L$ on $X$, we have a canonical isomorphism 
$$
L\cong L^\times\times_{\C^\times}\C_{\chi_1}.
$$
Therefore, a line bundle $L'$ is of the form $L_\chi$
for some $\C^\times$-character $\chi$ iff $L'$ is a power of $L$.
\end{lem}
\begin{proof}
Define 
$$
\rho:L^\times\times\C\ra L, ~~~(m,a)\mapsto am.
$$
Let $h\in H:=\C^\times$ acts by $h(m,1):=(mh^{-1},h)$ as usual.
Then $\rho:(mh^{-1},ha)\mapsto am$ for $h\in\C^\times$. It follows that $\rho$ descends to $L^\times\times_H\C_{\chi_1}\ra L$. It is clear that this is an isomorphism.

If $L'=lL$, then the isomorphism implies that $L'=L_{\chi_1^l}$. Conversely, if $L'=L_\chi$ with $\chi\in Hom(\C^\times,\C^\times)=\bra\chi_1\ket$, then $\chi=\chi_1^l$ and so $L'=lL$.
\end{proof}

\begin{thm} \label{obstruction}
(Obstruction to CY structure) Let $H-M{\br\pi\over\ra} X$ be a principal bundle. Then $M$ admits a CY structure iff there exists an $H$-character $\chi$ such that 
$$K_X\cong M\times_H\C_\chi.$$
In particular, if $L$ is line bundle on $X$, then $L^\times$ admits a CY structure iff $K_X$ is a power of $L$.
\end{thm}
\begin{proof}
By Theorem \ref{KX-theorem}, we have
$$
K_M/H\cong K_X\otimes L_{\chi_\gh}^{-1}.
$$ 
If $M$ has a CY structure $(\omega_M,\chi_M)$, then the left side is $L_{\chi_M}$, and so $K_X\cong L_{\chi_M\chi_\gh}$.

Conversely, suppose $K_X\cong L_\chi$ for some $H$-character $\chi$. Then by the isomorphism
$$
Pic_{G\times H}(M)\cong Pic_G(X),~~~F\mapsto F/H,~~~\pi^*E\leftarrow E,
$$
it follows that we have an $H$-equivariant isomorphism $K_M\cong M\times\C_{\chi_M}$ where $\chi_M:=\chi\chi_\gh^{-1}$, which implies that $K_M$ is trivial (but not necessarily $H$-equivariantly.) Moreover, the section $M\ra  M\times\C_{\chi_M}$, $m\mapsto(m,1)$, corresponds to a nowhere vanishing form $\omega_M\in H^0(K_M)^G$ which transforms by
$$
\Gamma_h\omega_M=\chi_M(h)\omega_M
$$
under $H$. So $M$ admits a CY structure.

Now specialize to the case $M=L^\times$, $H=\C^\times$. Then $M$ admits a CY structure iff there exists $\chi\in Hom(H,\C^\times)$, i.e. $\chi=\chi_1^l$ for some $l\in\Z$, such that 
$$K_X\cong M\times_H\C_\chi.$$
Then by the preceding lemma, the right side is $lL$, as desired.
\end{proof}

We now give a number of applications of the theorem. 

\begin{cor}\label{Calabi} 
(Calabi) $K_X^\times$ has a CY structure $(\hat\omega,\chi_1)$.
\end{cor}
\begin{proof}
This is the special case $L=K_X$ of the theorem.
\end{proof}

\begin{exm}
Projective space.
\end{exm}
\cut

Let $X=\PP^d$. By Theorem \ref{obstruction}, the principal $\C^\times$-bundle $\cO(k)^\times$ on $X$ admits a CY structure iff $k$ divides $d+1$.

\begin{exm}
A projective hypersurface as ambient space. 
\end{exm}
\cut

Let $X$ be a degree $n$ smooth hypersurface of $\PP^{d+1}$, $d\geq1$. Consider $L=\cO(k)$ on $X$. Then it follows immediately from Theorem \ref{obstruction} (and Lefschetz hyperplane) that $L^\times$ admits a CY structure iff $k$ divides $n-d-2$.

\begin{lem} 
Let $H-M\ra X$ be a principal bundle with $H=\C^\times$.  Then there exists a unique line bundle $L$ on $X$ such that $M\cong L^\times$, $H$-equivariantly.
\end{lem}
\begin{proof}
Set
$$
L=M\times_H\C_{\chi_1}.
$$
Then it is straightforward to show that $M\cong L^\times=M\times_H\C^\times$. The uniqueness part follows from Lemma \ref{L-star}.
\end{proof}

\begin{cor} (Classification of rank 1 CY bundles)
Put $H=\C^\times$. A principal bundle $H-M\ra X$ admits a CY structure iff $M$ is $H$-equivariantly isomorphic to $L^\times$ for some  line bundle $L$ such that $K_X$ is a power of $L$.
\end{cor}
\begin{proof}
This follows from the preceding lemma and Theorem \ref{obstruction}.
\end{proof}

We now prove the uniqueness of CY structures. 

\begin{lem} (Uniqueness lemma)
Let $(\omega,\chi)$ and $(\omega',\chi')$ be two CY structures on a given principal bundle $H-M\ra X$. If $\chi=\chi'$ then $\omega$ is a scalar multiple of $\omega'$. In particular, if the natural map $\Hom(H,\CC^\times)\ra Pic(X)$ is injective (cf. Proposition \ref{character-linebundle}) then up to scalar, there is at most one CY structure on $M$. 
\end{lem}
\begin{proof}
We have a (unique) function $f\in H^0(\cO_M^\times)$ such that
$$
\omega'=f\omega.
$$
For $h\in H$, we have $\Gamma_h\omega=\chi(h)\omega$ and $\Gamma_h\omega'=\chi'(h)\omega'$. This implies that
$$
\Gamma_h f=\chi^{-1}(h)\chi'(h) f.
$$
By Proposition \ref{function-section}, $f$ represents a nowhere vanishing global section of the line bundle $L_{\chi{\chi'}^{-1}}$ on $X$.  
Thus if $\chi=\chi'$, then $f$ defines a holomorphic function on $X$, hence must be constant, proving the first assertion.

Now suppose $\Hom(H,\CC^\times)\ra Pic(X)$ is injective. Then $\chi=\chi'$ by Theorem \ref{KX-theorem}, hence $\omega'$ is a scalar multiple of $\omega$ in this case as well.
\end{proof}

\begin{lem}
Let $(\omega,\chi)$ be a CY structure on $H-M\ra X$, $d\zeta$ be the standard coordinate 1-form on $\CC_{\chi\chi_\fh}$, and $x_1,..,x_q$ ($q=\dim H$) be independent vector fields generated by $H$ on $M\times\CC_{\chi\chi_\fh}$. Then the form
$$\hat\Omega:=\iota_{x_1}\cdots\iota_{x_q}(\omega\wedge d\zeta)$$ 
is $H$-basic, hence defines a nowhere vanishing top form on 
$$K_X\cong M\times_H\CC_{\chi\chi_\fh}.$$
\end{lem}
\begin{proof}
The form is obviously $H$-horizontal, and nowhere vanishing, since the vector fields $x_i$ are everywhere independent. We want to check that it is $H$-invariant. First note that $\Gamma_h\omega=\chi(h)\omega$ and $\Gamma_h d\zeta=\chi(h)^{-1}\chi_\fh(h)^{-1} d\zeta$. It follows that $\Gamma_h(\omega\wedge d\zeta)=\chi_\fh(h)^{-1}(\omega\wedge d\zeta)$. Since
\begin{eqnarray*}
\Gamma_h\iota_{x_1}\cdots\iota_{x_q}\Gamma_h^{-1}&=&\Gamma_h\iota_{x_1}\Gamma_h^{-1}\cdots\Gamma_h\iota_{x_q}\Gamma_h^{-1}\cr
&=&\iota_{Ad(h)x_1}\cdots\iota_{Ad(h)x_q}\cr
&=&\chi_\fh(h)\iota_{x_1}\cdots\iota_{x_q}\cr
\end{eqnarray*}
it follows that $\hat\Omega$ is $H$-invariant.
\end{proof}

\begin{cor}
The $\hat\Omega$ above, as a form on $K_X$, is a scalar multiple of a CY structure $\hat\omega$ of $K_X$.
\end{cor}
\begin{proof}
It suffices to show that the statement holds in $K_X^\times$. By the preceding lemma, $\hat\Omega$ is a nowhere vanishing top form on $K_X^\times$. Consider the $\C^\times$ action by the usual fiberwise scaling on $K_X$. Lifted to $M\times\C_{\chi\chi_\fh}$, it is just the usual scaling on $\C_{\chi\chi_\fh}$. Hence $d\zeta$ transforms by the $\C^\times$-character $\chi_1$. It follows that $\hat\Omega$ transforms by the character $\chi_1$ as well. Thus $(\hat\Omega,\chi_1)$ is a CY structure on the bundle $K_X^\times$. But $(\hat\omega,\chi_1)$ is also a CY structure on the same bundle, by Corollary \ref{Calabi}. By the Uniqueness Lemma, $\hat\Omega$ must be a scalar multiple of $\hat\omega$. 
\end{proof}

\begin{thm} (Uniqueness of CY structure) \label{uniqueness}
If a bundle $H-M\ra X$ admits a CY structure, it is unique up to scalar multiple. 
\end{thm}
\begin{proof}
Let $\omega_1$, $\omega_2$ be two CY structures on $M$. By the preceding corollary again, we can normalize them by a constant so that
\begin{equation}\label{e1}
\iota_{x_1}\cdots\iota_{x_q}(\omega_1\wedge d\zeta)=\iota_{x_1}\cdots\iota_{x_q}(\omega_2\wedge d\zeta)
\end{equation}
viewed as forms on $M\times\C$. Since ${\partial\over\partial\zeta}$ is a global vector field on $M\times\C$ (with no component on the $M$ side), we can apply $\iota_{\partial\over\partial\zeta}$ to eqn. \ref{e1} and get
\begin{equation}\label{e2}
\iota_{x_1}\cdots\iota_{x_q}\omega_1=\iota_{x_1}\cdots\iota_{x_q}\omega_2.
\end{equation}
Since multiplication by $\iota_{\partial\over\partial\zeta}$ kills the ${\partial\over\partial\zeta}$ component of each vector field $x_i$,
we can now regard the $x_i$ as the vector fields generated by the $H$-action on $M$ alone, and eqn.\ref{e2} is now an identity on $M$. Finally, since the $x_i$ are everywhere independent vector fields on $M$, \ref{e2} implies that
$$
\omega_1=\omega_2.
$$
%
\end{proof}

We now show that the CY bundle $\C^\times-K_X^\times\ra X$ is the universal one, in some sense.

\begin{prop} (Universal property)
Let $H-M{\br\pi_M\over\ra} X$ be any $H$-bundle with a CY structure $(\omega,\chi)$. Then the projection map $\pi_M$ factors through a canonical $H$-equivariant map $M\ra K_X^\times$, where $H$ acts on $K_X^\times$ fiberwise by the character $\chi$.
\end{prop}
\begin{proof}
The given CY structure determines a canonical isomorphism $K_X\cong M\times_H\C_{\chi\chi_\gh}$, by Theorem \ref{KX-theorem}. So, we may as well identify these two sides. Define $M\ra K_X^\times$, $m\mapsto[m,1]\in K_X$. For $h\in H$, we have
$$
mh^{-1}\mapsto[mh^{-1},1]=[m,\chi(h)^{-1}]=[m,1]\chi(h)^{-1}.
$$
This shows that our map is $H$-equivariant. Composing it with the projection $K_X^\times\ra X=M/H$, $[m,1]\mapsto[m]$, clearly yields $\pi_M:M\ra X$.
\end{proof}

\section{Poincar\'e residues for principal bundles}\label{sec-Poincare}

In this section, we continue to use the notations introduced earlier. We fix a CY bundle $H-M\ra X$, an $H$-character $\chi$, and put $L=L_\chi$ (cf. Proposition \ref{character-linebundle}.) 
Let $\sigma\in H^0(L)$ be any nonzero section. We now give an explicit construction of global meromorphic $d$-form ($d=\dim X$) with pole given by the divisor $\sigma=0$ (possibly singular) in $X$. From now on, given $x\in\fh$, we shall write the interior multiplication operator $\iota_{\xi_x}$ (acting on differential forms on $M$) simply as $\iota_x$.
Fix a nonzero element 
$$
x_1\wedge\cdots\wedge x_q\in\wedge^q \fh \blank(q=\dim H.)
$$
Put
\begin{eqnarray*}
\Omega&=&\iota_{x_1\wedge\cdots\wedge x_q}\omega_M=\iota_{x_1}\cdots\iota_{x_q}\omega_M\cr
\Omega_\sigma&:=&{1\over\sigma_M}\Omega\in H^0(K_M|(M-\{\sigma_M=0\})),
\end{eqnarray*}
 where $\sigma_M: M\ra\CC$ is the $\chi$-equivariant function representing $\sigma$. (Cf. Proposition \ref{function-section}.)

\begin{thm} (Poincar\'e residue)\label{residue}
For $\sigma\in H^0(L)-0$ and $\tau\in H^0(L+K_X)$, $\Omega_\sigma$ is a nowhere vanishing meromorphic $d$-form on $M-\{\sigma_M=0\}$ with pole along $\sigma_M=0$ such that
$$
\Gamma_g(\tau_M\Omega_\sigma)=(\Gamma_g\tau_M)\Omega_{g\circ\sigma\circ g^{-1}}\blank~(g\in G.)
$$
 Moreover $\tau_M\Omega_\sigma$ is $H$-basic, i.e.
$$
\iota_x(\tau_M\Omega_\sigma)=0,\blank\Gamma_h(\tau_M\Omega_\sigma)=\tau_M\Omega_\sigma\blank(x\in\fh,~~ h\in H.)
$$
Therefore, $\tau_M\Omega_\sigma$ defines a meromorphic $d$-form on $X-\{\sigma=0\}$ with pole given by the divisor $\sigma=0$.
\end{thm}
\begin{proof}
Since $\xi_{x_1},..,\xi_{x_q}$ are pointwise independent vector fields on $M$, and $\omega_M$ is a nowhere vanishing holomorphic top form on $M$, it follows that $\Omega$ above is a nowhere vanishing form of degree $d=\dim M-q$. Since the vector fields $\xi_{x_i}$ are $G$-invariant (because the $G$- and the $H$-actions on $M$ commute), it follows that $\Gamma_g\Omega=\Omega$. Since $L$ is a $G$-equivariant bundle, there is an induced linear action on $H^0(L)$. By Proposition \ref{function-section}, our first assertion follows.

Next, we show the $H$-basic property. By construction of $\Omega$, it is clear that $\iota_x\Omega=0$ for all $x\in\fh$. Since $\iota_x$ is a derivation acting trivially on functions, It follows $\iota_x(\tau_M\Omega_\sigma)=0$ for all $x\in\fh$. To see that $\tau_M\Omega_\sigma$ is $H$-invariant, let $h\in H$. First note that
\begin{eqnarray*}
\Gamma_h\Omega&=&\Gamma_h\iota_{x_1}\Gamma_h^{-1}\cdots\Gamma_h\iota_{x_q}\Gamma_h^{-1}\Gamma_h\omega_M\cr
&=&\iota_{Ad(h)x_1}\cdots\iota_{Ad(h)x_q}\chi_M(h)\omega_M\cr
&=&\chi_\fh(h)\chi_M(h)\iota_{x_1}\cdots\iota_{x_q}\omega_M\cr
&=&\chi_\fh(h)\chi_M(h)\Omega.
\end{eqnarray*}
On the other hand, by Theorem \ref{KX-theorem}, we have $\tau\in H^0(L+K_X)=H^0(L_{\chi\chi_M\chi_\fh})$, hence
$$
\tau_M(mh^{-1})=\chi(h)\chi_M(h)\chi_\fh(h)\tau_M(m).
$$
It follows that
$$
(\Gamma_h{\tau_M\over\sigma_M})(m)={\tau_M(mh)\over\sigma_M(mh)}=\chi_M(h)^{-1}\chi_\fh(h)^{-1}{\tau_M(m)\over\sigma_M(m)}.
$$
This implies that $\Gamma_h(\tau_M\Omega_\sigma)=\tau_M\Omega_\sigma$.

The last assertion of the theorem follows from the $H$-basic property of $\tau_M\Omega_\sigma$.
\end{proof}

\begin{rem}
One important way in which the preceding theorem will be applied is when $\tau\in H^0(L+K_X)$ transforms under a group $G\subset\Aut X$ by a given $G$-character, in which case $G$ ``moves'' the form $\tau_M\Omega_\sigma$ essentially by varying only its pole, according to the theorem.
\end{rem}

\section{Examples: old and new}\label{sec-examples}

This section is intended to offer some new perspectives on a number of  known constructions of Poincar\'e residues and period integrals. By interpreting them in the context of the earlier sections, we can now unify all of them on more or less equal footing.

\begin{exm}
$M=\CC^{d+1}-0$, $H=\CC^\times$, $G=SL_{d+1}$, $K=G\times H$, $X=M/H=\PP^d$.
\end{exm}
\cut
We will regard $M$ as a set of column vectors so that $G$ acts on the left as usual, and $H$ acts on the right by scaling $h:M\ra M$, $m\mapsto mh^{-1}$. From this we have $dh({\partial\over\partial m_i})=h^{-1}{\partial\over\partial m_i}$ and $\delta h^{-1}(dm_i)=h dm_i$. Let $\omega_M=dm_0\wedge\cdots\wedge dm_d$, which is clearly $G$-invariant, and 
$$
\Gamma_h\omega_M=h^{d+1}\omega_M.
$$ 
So, the $H$-character for $\CC\omega_M$ is  $\chi_M(h)=h^{d+1}$. The vector field generated by $H$ on $M$ is $x=-\sum m_i{\partial\over\partial m_i}$, and so
$$
\Omega=\iota_x\omega_M=-\sum (-1)^i m_i dm_0\wedge\cdots\widehat{dm_i}\cdots\wedge dm_d.
$$
The $H$-character for $\fh$ is trivial since $H$ is abelian. By Theorem \ref{KX-theorem}, it follows that
$$
K_X\cong L_{\chi_M}.
$$
Note that the coordinate functions $m_i$ on $M$ represents the sections of $\cO(1)$ of $X$. Since
$$
(mh^{-1})_i=h^{-1}m_i
$$
it follows that the $H$-character for the line bundle $\cO(1)$ is $\chi(h)=h^{-1}$, hence $\chi_M=\chi^{-d-1}$ and $K_X=\cO(-d-1)$, as expected.

\bs
\begin{exm} \label{flagbundle-example}
Partial (Type A) flag variety $F(d_1,...,d_r,n)$ revisited.
\end{exm}
\cut
Put
\begin{eqnarray*}
M&=&M_{d_{r+1},d_r}\times\cdots\times M_{d_2,d_1} \blank(n:=d_{r+1}>\cdots>d_1>d_0:=0)\cr
H&=&GL_{d_r}\times\cdots\times GL_{d_1}\cr
G&=&SL_n\cr
X&=&M/H=F(d_1,..,d_r,n).
\end{eqnarray*}
Here $M_{a,b}$ is the space of $a\times b$ matrices of full rank. Here $G$ and $H$ act on $M$ by
\begin{eqnarray*}
&&G\times M\ra M,~~~(g,m_r,..,m_1)=(gm_r,..,m_1)\cr
&&H\times M\ra M,~~~(h_r,..,h_1,m_r,..,m_1)\mapsto(m_rh_r^{-1},h_rm_{r-1}h_{r-1}^{-1},..,h_2m_1h_1^{-1}).
\end{eqnarray*}
Put $M_i=M_{d_{i+1},d_i}$ and let $\omega_i$ be the coordinate top form on $M_i$ (with respective to some chosen order, say, the lexicographic order, on the coordinates of $M_i$). Put
$$
\omega_M=\omega_r\wedge\cdots\wedge\omega_1.
$$
Then $\omega_M$ is clearly $G$-invariant. For $h=(h_r,..,h_1)\in H$, we have
$$
\Gamma_h\omega_M=\chi_M(h)\omega_M
$$
where $\chi_M(h):=\prod_{i=1}^r\Det_{d_i}(h_i)^{d_{i+1}-d_{i-1}}$. Here $\det_d:GL_d\ra\CC^\times$ is the determinant function. Since the Lie algebra of $GL_d$ is $\fs\fl_d\oplus\CC$, it follows that
$\wedge^q\fh$ ($q=\dim H$) is a trivial representation of $H$, i.e. $\chi_\fh=1$ in this case. By Theorem \ref{KX-theorem}, it follows that
$$
K_X\cong L_{\chi_M}.
$$

\begin{prop}\label{typeA-Picard} (Cf. Proposition \ref{Popov}) For the principal $H$-bundle $H-M\ra X$ given in the preceding example, the map
$$\Hom(H,\CC^\times)\ra Pic_G(X),~~~\chi\mapsto [L_\chi]$$
is isomorphism.
\end{prop}
\begin{proof}
The group $G=SL_n$ is simple, connected, simply-connected, and $X$ is compact and $G$-homogeneous. It follows the group $Pic_G(X)\cong Pic(X)$ is free abelian generated by the $G$-equivariant line bundles corresponding to some fundamental weights of $G$ \cite{bott}\cite{Popov}. They are given by the highest weights $\lambda_{d_i}$ of the representations $\wedge^{d_i}\CC^n$, $1\leq i\leq r$. In fact, the line bundle corresponding to $\lambda_{d_i}$ is precisely the pullback of the hyperplane bundle $\cO_i(1)$ on the Grassmannian $G(d_i,n)$, under the canonical projection $X\ra G(d_i,n)$.  In turn, $\cO_i(1)$ is the pullback of the universal hyperplane bundle, under the Pl\"ucker embedding $G(d_i,n)\into \PP(\wedge^{d_i}\CC^n)$. We shall henceforth identify $Pic_G(X)$ with the lattice generated by the $\lambda_{d_i}$. On the other hand, the character group $\Hom(H,\CC^\times)$ is free abelian generated by the characters $\chi_i:=\det_{d_i}^{-1}$ (viewed as functions on $H$.) So, to complete the proof, it suffices to show that the line bundle $L_{\chi_i}$ is isomorphic to the pullback of $\cO_i(1)$.

The projection $X\ra G(d_i,n)$ can be realized as the 
$$
M/H\ra M_{n,d_i}/GL_{d_i},~~~[m_r,...,m_1]\mapsto [m_r\cdots m_i].
$$ 
Fix an index set $J=(1\leq j_1<\cdots<j_{d_i}\leq n)$ and let $m_J$ to be the $d_i\times d_i$ minor of $m\in M_{n,d_i}$ consisting of the rows indexed by $J$, and consider the Pl\"ucker coordinate for $G(d_i,n)$
$$
\sigma_J:M_{n,d_i}\ra\CC,~~m\mapsto \det(m_J).
$$
For $h_i\in GL_{d_i}$, we have $\sigma_J(mh_i^{-1})=\det_{d_i}(h_i)^{-1}\sigma_J(m)$. By Proposition \ref{function-section}, $\sigma_J$ represents a (nonzero) section of the line bundle $L_{\det_{d_i}^{-1}}$ on $G(d_i,n)=M_{n,d_i}/GL_{d_i}$. It follows that $L_{\det_{d_i}^{-1}}=\cO_i(1)$. But we have the commutative diagram of bundles
$$
\begin{array}[c]{ccc} 
M\times_H\CC_{\chi_i}&\ra&M_{n,d_i}\times_{GL_{d_i}} L_{\det_{d_i}^{-1}}\cr
 \downarrow& &\downarrow\cr
 X&\ra&G(d_i,n)
\end{array}
$$
which means that $L_{\chi_i}$ is isomorphic to the pullback of $L_{\det_{d_i}^{-1}}=\cO_i(1)$, via the projection $X\ra G(d_i,n)$.
\end{proof}

In \cite{LSY1}, we have also carried out an explicit construction of the period integrals for CY complete intersections in the partial flag variety $X=F(d_1,..,d_r,n)$. In a later section, we will generalize the construction to an arbitrary homogeneous space, and for CY as well as  general type complete intersections.

\bs
\begin{exm} Orthogonal (Types B and D) and isotropic (Type C) flag varieties $OF(d_1,...,d_r,n)$, $IF(d_1,..d_r,n)$.
\end{exm}
\cut
Let $M, G,H,X$ as in the preceding example, and assume that $d_r\leq n/2$. Put
\begin{eqnarray*}
M_1&=&\{(m_r,..,m_1)\in M|m_r^tJ_1m_r=0\}\cr
G_1&=&SO_n\cr
X_1&=&M_1/H=OF(d_1,..,d_r,n).
\end{eqnarray*}
Here $SO_n$ is the orthogonal group: the group of isometry of the symmetric bilinear form $(u,v)\mapsto u\cdot J_1 v$ on $\CC^n$, where $J_1$ is the matrix (in column form) $J_1=[e_n,..,e_1]$. Here $e_1,..,e_n$ is the standard basis of $\CC^n$. It is easy to show that $M_1$ is an affine complete intersection of codimension $c=\half d_r(d_r+1)$ in $M$. In other words, the components of the quadratic equations $m_r^tJm_r=0$ above define transversal affine hypersurfaces in $M$. 

Next, we put
\begin{eqnarray*}
M_2&=&\{(m_r,..,m_1)\in M|m_r^tJ_2m_r=0\}\cr
G_2&=&Sp_{2l}~~(n=2l)\cr
X_2&=&M_2/H=IF(d_1,..,d_r,n).
\end{eqnarray*}
Here $Sp_{2l}$ is the symplectic group: the group of isometry of the skew symmetric bilinear form $(u,v)\mapsto u\cdot J_2 v$ on $\CC^{2l}$, where $J_2$ is the matrix (in column form) $J_1=[e_{2l},..,e_l,-e_{l-1},..,-e_1]$. Again, $M_2$ is an affine complete intersection of codimension $c=\half d_r(d_r-1)$ in $M$. 

We can apply Theorems \ref{KX-theorem} and \ref{residue} to the principal $H$-subbundles $M_1\subset M$ above to compute $K_{X_1}$ and Poincar\'e residues for CY hypersurfaces in $X_1$. One way is to first construct a holomorphic top form $\omega_{M_1}\in H^0(K_{M_1})$ by taking Poincar\'e residue of $\omega_M$ over the affine complete intersection $m_r^tJ_1m_r=0$. The same can be said about $X_2$.
However, we will not carry that out here. Instead, later we will apply Theorems \ref{KX-theorem} and \ref{residue} to compute $K_X$ and Poincare residues in an entirely different way. In fact, we will do so for all homogeneous spaces $X$ of any complex Lie group $G$.

\bs
\bs
\begin{exm}\label{toric-example}
Toric manifolds revisited.
\end{exm}
\cut
Let $X$ be a toric manifold, associated with a fan $\Sigma$ in $\ZZ^n$. We briefly recall a construction of $X$ as a quotient \cite{Audin}\cite{Cox}, and a description of its line bundles. Let $\nu_1,..,\nu_t\in\ZZ^n$ be the integral generators of the 1-cones in $\Sigma$. For simplicity, we shall assume that {\it no hyperplane in $\ZZ^n$ contains all $t$ points $\nu_i$, and that every maximal cone in $\Sigma$ is $n$-dimensional.} Any complete fan, for example, satisfies this assumption. Put
\begin{eqnarray*}
\tilde G&=&(\CC^\times)^t\cr
G&=&\{g=(g_1,..,g_t)\in \tilde G|g_1\cdots g_t=1\}\cr
H&=&\{g=(g_1,..,g_t)\in \tilde G|\prod_{j=1}^t g_j^{\bra \mu,\nu_j\ket}=1~~~\forall \mu\in(\ZZ^n)^*\}\cr
T&=&(\CC^\times)^n\cr
\cL&=&\{l=(l_1,..,l_t)\in\ZZ^t|\sum_j l_j\nu_j=0\}.
\end{eqnarray*}
We have exact sequences
\begin{eqnarray*}
& &H\into\tilde G\onto T\cr
& &h\mapsto h,~~~g\mapsto (\prod_j g_j^{\bra e_1^*,\nu_j\ket},..,\prod_j g_j^{\bra e_n^*,\nu_j\ket})\cr
& &(\ZZ^n)^*\into(\ZZ^t)^*\onto\cL^*\cr
& &\mu\mapsto\sum_j\bra \mu,\nu_j\ket D_j,~~~\sum_j a_j D_j\mapsto (\lambda_a:l\mapsto \sum_j l_j a_j).
\end{eqnarray*}
Here $e_1^*,..,e_n^*$ form the basis of $(\ZZ^n)^*$ dual to the standard basis of $\ZZ^n$, and $D_1,..,D_t$ the basis of $(\ZZ^t)^*$ dual to the standard basis $D_1^*,..,D_t^*$ of $\ZZ^t$. Under the assumption above, the restriction of $\tilde G\onto T$ to $G$ is also onto. 
\comment{See subsection ``Why replace $\tilde G = (\C^\times)^p$ by $G$: $\omega_M$ is not $\tilde G$-invariant''.}

Let $G\times H$ act on $\CC^t$ by
$$
G\times H\times\CC^t\ra\CC^t,~~(g,h,m)\mapsto gmh^{-1}:=(g_1m_1h_1^{-1},..,g_tm_th_t^{-1}).
$$
 Obviously, $H$ acts freely and properly on $(\CC^\times)^t$. In particular, there is a unique maximal open subset $M$ of $\CC^t$ containing $(\CC^\times)^t$ such that the action on $M$ is free and proper. It can be described as follows. If $\sigma=\sum_{i\in I}\RR_+\nu_i\in\Sigma$, we put $\hat\sigma=\sum_{i\in I}\RR_+D_i^*\subset\RR^t$. Let $\hat\Sigma$ be the collection of all such $\hat\sigma$. This is a fan in $\RR^t$ whose corresponding toric variety is $M$. There is an induced action of $G$ on $M/H=M/(1,H)$. Since $G$ is abelian, $H\cap G$ as a subgroup of $G$ obviously act trivially on $M/H$, hence the $G$-action on $M/H$ descends to a $T\cong G/(H\cap G)$-action.

\begin{thm}\label{toric-Picard}
\cite{Fulton}\cite{Oda}
$$\Pic(X)\cong H^2(X,\ZZ)\cong (\ZZ^t)^*/(\ZZ^n)^*\equiv\cL^*.$$
\end{thm}

\begin{thm}
\cite{Audin}\cite{Cox} There is a $T=G/(H\cap G)$-equivariant isomorphism
$$M/H\cong X.$$
\end{thm}

Let $z_1,..,z_t$ be the standard coordinates of $\CC^t$. Then $M$ is of the form $\CC^t-Z$ where $Z$ is a certain subvariety of $ \cup_j\{z_j=0\}\subset\CC^t$. In fact, $Z$ can be described explicitly in terms of the fan $\Sigma$ \cite{Cox}. Moreover, the cohomology class $[D_j]\in H^2(X,\ZZ)$ is  the first Chern class of the line bundle corresponding to the divisor $H\cdot\{z_j=0\}$ in $M$.

By shifting the focus from $X$ to $M$, viewed as a $G\times H$-space and a $G$-equivariant principal $H$-bundle over $X$, we can apply results of  \ts\ts\ref{sec-adjunction}-\ref{sec-Poincare} to $M$ as a special case. This will allow us to reinterpret $\Pic(X)$ in terms of $H$-characters, and to carry out the Poincar\'e residue construction for $M$. The latter will recover a well-known construction of periods of CY hypersurfaces in $X$, and the differential equations governing them, when $X$ is compact. But the construction makes sense even without the compactness assumption.

\begin{prop} (Cf. Propositions \ref{function-section}, \ref{Popov}, \ref{typeA-Picard})
The map
$$\Hom(H,\CC^\times)\ra Pic(X),~~~\chi\mapsto [L_\chi]$$
is isomorphism.
\end{prop}
\begin{proof} \comment{See Subsection ``From $H$-characters to line bundles'' for more details.}
With the identification $\cL_\CC/\cL\equiv H$, $\bar l=l+\cL\equiv (e^{2\pi i l_1},..,e^{2\pi i l_t})$, we have
$$
\cL^*{\br\sim\over\ra} \Hom(H,\CC^\times)\equiv \Hom(\cL_\CC/\cL,\CC^\times),~~~\lambda\mapsto(\chi_\lambda:\bar l\mapsto e^{2\pi i\lambda(l)}).
$$
Composing this with $\Pic(X){\br\sim\over\ra}\cL^*$, $\sum_j a_j[D_j]\mapsto(\lambda_a:l\mapsto\sum_j l_j a_j)$, Theorem \ref{toric-Picard}, we get
$$
\Pic(X){\br\sim\over\ra}\Hom(H,\CC^\times),~~~\sum_i a_j[D_j]\mapsto\chi_{\lambda_a}.
$$
To complete the proof, it suffices to show that
\begin{equation}\label{to-show}
[L_{\chi_{\lambda_a}}]=-\sum_i a_j[D_j].
\end{equation}
It is enough to show this for a set of generators, say $[L_{\chi_j}]\cong [D_j]$  where 
$$\chi_j(\bar l)=e^{-2\pi i l_j}\blank (j=1,..,t).$$

Since the function $z_j:M\ra\CC$ satisfies $z_j(m{\bar l}^{-1})=e^{-2\pi i l_j}m_j=\chi_j(\bar l)z_j(m)$, it defines a section of $L_{\chi_j}$ by Proposition \ref{function-section}. Hence its zero locus is the divisor representing $L_{\chi_j}$. But the zero locus is exactly the $H$ orbit of $z_j=0$ in $M$. This divisor represents the class $[D_j]$.
\end{proof}

Next, we apply results of \ts\ts\ref{sec-adjunction}-\ref{sec-Poincare} to construct Poincar\'e residues for CY hypersurfaces in $X$. This example has been previously considered \cite{Danilov}\cite{BatyrevCox1994}.
 Put
$$
\omega_M=dz_1\wedge\cdots\wedge dz_t,~~~\chi_M(h)=h_1\cdots h_t\blank (h\in H.)
$$
Then $\omega_M$ is $G$-invariant and $\Gamma_{h}\omega_M=\chi_M(h)\omega_M$. Since $H$ is abelian, $K_X\cong L_{\chi_M}$  by Theorem \ref{KX-theorem}. 
Since $\chi_M=(\chi_1\cdots\chi_t)^{-1}$ and $[L_{\chi_j}]= [D_j]$,  it follows that $[K_X]= -\sum_j[D_j]$, as expected.

Fix a $\ZZ$-basis $l^1,..,l^q$ ($q=t-n$) of the lattice $\cL$. The vector field generated on $M$ by $l^i\in\cL_\CC$, the Lie algebra of $H\equiv\cL_\CC/\cL$, is
$$
x_i=\sum_j l^i_j z_j {\partial \over\partial z_j}.
$$

\begin{prop} $\Omega=\iota_{x_1}\cdots\iota_{x_q}\omega_M$ in Theorem \ref{residue} is a nonzero constant multiple of the holomorphic $n$-form on $M$
$$
\Omega'=z_1\cdots z_t\wedge_{k=1}^n\sum_{j=1}^t\bra e_k^*,\nu_j\ket {dz_j\over z_j}.
$$
\end{prop}
\begin{proof}
The apparent poles of $\Omega'$ are clearly removable. Put $\omega_{q+k}=\sum_{j=1}^t\bra e_k^*,\nu_j\ket {dz_j\over z_j}$ ($k=1,..,n$.) We'll show that for some constant $c$, we have on the dense subset $\tilde G=(\CC^\times)^t\subset M$, 
$$
{\Omega\over z_1\cdots z_t}=c\omega_{q+1}\wedge\cdots\wedge\omega_t.
$$

Since the $\nu_j$ span $\ZZ^n$, $\omega_{q+1},..,\omega_t$ are pointwise independent. Since $l^i\in\cL$, it follows that
$$
\iota_{x_i}\omega_{q+k}=\sum_j l_j^i\bra e_k^*,\nu_j\ket=0 \blank(1\leq i\leq q,~1\leq k\leq n.)
$$
Since $x_1,..,x_q$ are $\tilde G$-invariant vector fields generated by the $H$-action, and since the $\omega_{q+k}$ are $\tilde G$-invariant forms on $\tilde G$, we can find $\tilde G$-invariant vector fields $x_{q+1},..,x_t$ and 1-forms $\omega_1,..,\omega_q$ such that, for some $c$,
\begin{eqnarray*}
\iota_{x_i}\omega_k&=&\delta_{ik}\blank(1\leq i,k\leq t)\cr
{dz_1\wedge\cdots\wedge dz_t\over z_1\cdots z_t}&=&c\omega_1\wedge\cdots\wedge\omega_t.
\end{eqnarray*}
Applying $\iota_{x_1}\cdots\iota_{x_q}$ to the last equation completes the proof.
\end{proof}

By Theorem \ref{residue},
$$
\Gamma_g\Omega_\sigma=\Omega_{g\circ\sigma\circ g^{-1}}.
$$
Next, let us restrict $\Omega_\sigma$ to the torus $\tilde G\subset M$, and examine it more closely. Let $y_1,..,y_n$ be the standard coordinate functions on $T=(\CC^\times)^n$. We can pull them back to $\tilde G$ via the map $\tilde G\onto T$, $g\mapsto (\prod_j g_j^{\bra e_1^*,\nu_j\ket},..,\prod_j g_j^{\bra e_n^*,\nu_j\ket})$, express the pullback (still denoted by $y_k$) in terms of the functions $z_j$ on $\tilde G$, and get
$$
y_k=\prod_j z_j^{\bra e_k^*,\nu_j\ket}.
$$
It follows that ${dy_k\over y_k}=\omega_{q+k}$ and by the preceding proposition,
$$
{\Omega\over z_1\cdots z_t}=c{dy_1\over y_1}\wedge\cdots\wedge{dy_n\over y_n}.
$$
Descending to $X$, we can view this as a form defined on the open dense $T$-orbit $T_0\subset X$. Consider those sections $\sigma$ of $-K_X$ which can be represented as a polynomial in the $z_1,..,z_t$. (Note that if $X$ is assumed compact, then every section can be represented as a polynomial in the $z_i$.)  A monomial $z^\alpha=z_1^{a_1}\cdots z_t^{a_t}$ is a section iff $\sum a_j[D_j]=-[K_X]=\sum_j[D_j]$. \comment{See Subsection ``Sections of a line bundle''.}
This holds iff there is a unique $\mu\in (\ZZ^n)^*$ such that $a_j-1=\bra \mu,\nu_j\ket$ for all $j$. \comment{This follows from the cohomology ring relation of $X$.}  In this case,
$$
z^\alpha=z_1\cdots z_t y^\mu.
$$ 
Thus $\sigma$, restricted to $T_0$, is $z_1\cdots z_t$ times a Laurent polynomial $f_\sigma$ in the $y_k$.
For $X$ a Fano toric variety, it is known \cite{Batyrev} that certain periods of the CY hypersurface $\sigma=0$ can be expressed as integral of the meromorphic form
$$
{1\over f_\sigma(y)}{dy_1\over y_1}\wedge\cdots\wedge{dy_n\over y_n}
$$
over cycles in $T_0$. We see that this form agrees with the meromorphic form $\Omega_\sigma|T_0$ in Theorem \ref{residue}, in this special case.

\bs
Next, we show that in a Fano manifold $X^d$, essentially all holomorphic top forms of CY (and general type) complete intersections can be obtained by means of Poincar\'e residues. First we have the following easy extension of the Poincar\'e residue map \cite{GH} for hypersurfaces to complete intersections.

\begin{lem} (Codimension $s$ Poincar\'e residue sequence) Let $Y_1,..,Y_s$ be smooth hypersurfaces in $X$ such that $Y=Y_1\cap\cdots\cap Y_s$ is smooth of codimension $s$ in $X=X^d$. Then we have the following exact sequence of sheaves:
$$
\sum_i\Omega^d_X(\cup_{j\neq i} Y_j)\into\Omega^d_X(\cup_i Y_i)\onto\Omega_Y^{d-s}
$$
where the first map is the inclusion map and the second map is given by
$$
{g\over z_1\cdots z_s}dz_1\wedge\cdots\wedge dz_d\mapsto g dz_{s+1}\wedge\cdots\wedge dz_d|Y,
$$
where the $z_i$ are local coordinates on $X$ centered at $p\in Y$, where a local equation of $Y_i$ is $z_i=0$.
\end{lem}

\begin{prop} \label{PR-isomorphism}
Suppose $X$ is Fano and  $Y_1,..,Y_s$ are hypersurfaces in $X$ such that for any $J\subsetneq\{1,..,s\}$, $\cap_{j\in J}Y_j$ is also Fano and smooth of codimension $|J|$ in $X$. Then the Poincar\'e residue map induces an isomorphism
$$
H^k(\Omega^d_X(\cup_i Y_i))\cong H^k(\Omega_{\cap_i Y_i}^{d-s})
$$
for $k<d-s$.
\end{prop}
\begin{proof}
We proceed by induction on $s$. For $s=1$, this follows from the codimension 1 Poincar\'e residue sequence and the fact that $H^0(\Omega^d_X)=0$ because $X$ is Fano. For simplicity, we write $\Omega\equiv\Omega^d_X$ and $W_i\equiv\cup_{j\neq i}Y_j$ in the following discussion. For $s>1$, the Koszul complex gives a resolution of 
$$A_0:=\sum_i\Omega(W_i)$$ 
whose $p$th term ($1\leq p\leq s$) is
$$
A_p:=\oplus_{|I|=p}\Omega(W_{i_1})\cap\cdots\cap\Omega(W_{i_p})=\oplus_{|J|=s-p}\Omega(\cup_{j\in J} Y_j)
$$
where the first sum is over all induces $I=(i_1<\cdots<i_p)$, and the second sum over the $J=\{1,..,s\}-I$. By induction on $s$, we have
$$
H^k(\Omega(\cup_{j\in J}Y_j))=H^k(\Omega^{d-s+p}_{Y_J})
$$
for $k<d-s+p$, where $Y_J:=\cap_{j\in J}Y_j$. Since each $Y_J$ is Fano, it follows that the right side is zero, and this yields
$$
H^k(A_p)=0,\blank k<d-s+p.
$$
Since $A_\bullet\ra A_0$ is a resolution, the usual long exact sequence argument shows that
$$
H^k(A_0)=0,\blank k<d-s+1.
$$
But $A_0$ is the first term of the codimension $s$ Poincar\'e residue sequence. Now our assertion follows from the preceding lemma.
This completes our proof.
\end{proof}

\section{A general global Poincar\'e residue formula}\label{sec-globalPR}

For the rest of this paper, unless stipulated otherwise, 

{\it $X$ will be a $d$-dimensional compact complex manifold.} 

The results in \ts\ref{sec-CYstructures} settle the existence and uniqueness questions for CY structures on a given principal $H$-bundle, up to scalar multiple. They do not, however, fix the normalization of the CY structure. In this section, we begin with a local construction of a CY structure on $K_X$, due to Calabi \cite{Calabi} (cf. Corollary \ref{Calabi}), that does determine a {\it canonical} normalization. Remarkably, this normalization turns out to be exactly the one that makes the global Poincar\'e residue formula (Theorems \ref{residue}, \ref{globalPR}) agrees, on the nose, with the map $\R$ in Definition \ref{period-integrals}. In addition, the local construction will allow us to compare global Poincar\'e residue formulas derived from different CY bundles. This was necessary to prove that CY structures on any principal bundle over $X$ are invariant under $\Aut X$ (Theorem \ref{uniqueness}.)

\begin{prop} (Calabi) $K_X$ admits a nowhere vanishing holomorphic top form $\hat\omega$.
\end{prop}
\begin{proof}
Let $\pi:K_X\ra X$ be the projection. Let $w=(w_1,..,w_d):U\ra\CC^d$ be a coordinate chart on $X$. This induces a chart (local trivialization)
$$
(z_w,w):\pi^{-1}(U)\ra\CC\times\CC^d,~~\eta\mapsto(\eta({\partial\over\partial w_1},..,{\partial\over\partial w_d}),w(\pi\eta))
$$
on $K_X$. Define the nowhere vanishing top form on $K_
X|U$:
$$
\hat\omega=dz_w\wedge dw_1\wedge\cdots\wedge dw_d.
$$
Clearly, under a coordinate change from $w$ to $w'$, the function $z_w$ and the form $dw_1\wedge\cdots\wedge dw_d$ transform respectively by the determinant of Jacobian matrices of $w'\circ w^{-1}$ and its inverse, which cancel out each other to leave $\hat\omega$ unchanged. It follows that $\hat\omega$ is independent of the choice of coordinates $w$, implying that $\hat\omega$ is globally defined.
\end{proof}

Calabi showed that if $X$ is compact K\"ahler, then $K_X$ admits a Ricci flat complete K\"ahler metric.

\begin{defn} 
We shall call $\hat\omega$ above the canonical form of $K_X$.
\end{defn}

We shall also denote by $\pi:K_X^\times\ra X$ the projection map. We now apply the results in  \ts\ts \ref{sec-adjunction}-\ref{sec-Poincare} to this particular principal bundle $\CC^\times-K_X^\times\ra X$, with $G$ being any complex Lie subgroup of the automorphism group $\Aut X$. Note that by the convention in \ts\ref{sec-adjunction}, $h\in H$ acts by $h\cdot m=m h^{-1}$. 

\begin{thm}\label{Calabi-form}
The canonical form $\hat\omega$ of $K_X$ is invariant under $G=\Aut X$. Moreover, for $h\in H=\CC^\times$, we have
$$
\Gamma_h\hat\omega=\chi_1(h)\hat\omega
$$
where $\chi_1(h)=h$. Therefore, $(\hat\omega,\chi_1)$ is a $G$-equivariant CY structure on $K_X^\times$. In particular, $\iota_{x_0}\hat\omega$ defines a $G$-invariant $\CC^\times$-horizontal $d$-form on $K_X^\times$, where $x_0$ is the vector field generated by $-1\in\CC=Lie(H)$. Finally, in local coordinates $w$ we have
$$
\iota_{x_0}\hat\omega=z_wdw_1\wedge\cdots\wedge dw_d.
$$
\end{thm}
\begin{proof}
Put $H=\CC^\times$ and $M=K_X^\times$. Let $w:U\ra \CC^d$ be a chart on $X$, and $(z_w,w)$ the induced chart on $M$, as before.
Recall that
$$
\hat\omega|U=dz_w\wedge dw_1\wedge\cdots\wedge dw_d.
$$
Since $z_w$ is a linear coordinate on the fiber on which $H$ acts by scaling, and since $\hat\omega$ is linear in $z_w$,  it is clear that $\Gamma_h\hat\omega=h\hat\omega$. The last assertion follows from that $x_0=z_w{\partial\over\partial z_w}$.

It remains to verify that $\hat\omega$ is $G$-invariant. Let $g\in \Aut X$; we will compute $\Gamma_g\hat\omega$. 
Since $g$ is holomorphic, $w^g:=w\circ g^{-1}:gU\ra\CC^d$ is a chart on $X$. So, it induces a chart on $M$:
$$
(z_{w^g},w^g):\pi^{-1}(gU)\ra\CC\times\CC^d.
$$
On the other hand $g\in\Aut X$ induces a holomorphic map $g:M\ra M$, hence a chart
$$
(z_w^g,w^g):\pi^{-1}(gU)\ra\CC\times\CC^d
$$
where $z_w^g=z_w\circ g^{-1}$. We will show that this chart is equal to $(z_{w^g},w^g)$, i.e. 
$$
(*)\blank z_w^g=z_{w^g}.
$$
Note that this completes the proof. For then we have
$$
(g^{-1})^*\hat\omega|U=(g^{-1})^*(dz_w\wedge dw_1\wedge\cdots\wedge dw_d)
=dz_w^g\wedge dw_1^g\wedge\cdots\wedge dw_d^g.
$$
Since the first factor of this form is $dz_{w^g}$ by (*), this form is equal to $\hat\omega|gU$. It then follows that
$$
(g^{-1})^*\hat\omega=\hat\omega
$$
as desired. 

We now return to (*). By definition, for $\eta'\in \pi^{-1}(gU)=K_X^\times|gU$,
$$
z_{w^g}(\eta')=\eta'({\partial\over\partial w_1^g},..,{\partial\over\partial w_d^g}),\blank
z_w^g(\eta')=z_w(g^{-1}_M\eta')=(g^{-1}_M\eta')({\partial\over\partial w_1},..,{\partial\over\partial w_d}).
$$
Since $g^{-1}_M\eta'=\delta g(\eta')=\eta'\circ dg$, proving (*) reduces to checking that
$$
dg({\partial\over\partial w_i})={\partial\over\partial w_i^g}
$$
on $X$, which is obviously true. This completes the proof.
\end{proof}


We now return to the set-up in the Introduction, spell out the assumptions, and fix notations that will be repeatedly used throughout the rest of the paper.

\begin{itemize}
\item{\bf X1.} Let $X$ be a compact complex manifold of dimension $d$, and fix a complex Lie subgroup $G\subset\Aut X$ (possibly trivial). Let $L_1,..,L_s$ ($s<d$) be nontrivial $G$-equivariant line bundles on $X$ such that $H^0(L+K_X)\neq0$, that each of the linear systems $V_i^*:=H^0(L_i)$ is base point free, 
and that there is an $H$-character $\chi$ (cf. Proposition \ref{character-linebundle}) such that
$$
L:=L_1+\cdots+L_s=L_\chi.
$$
\item{\bf X2.} Assume that for general sections $\sigma_i\in V_i^*$, the intersection $Y_{\sigma_1,...,\sigma_s}:=Y_{\sigma_1}\cap\cdots\cap Y_{\sigma_s}$ ($Y_{\sigma_i}=\{\sigma_i=0\}$) is smooth of codimension $s$. Put 
$$V:=V_1\times\cdots\times V_s,\blank\blank~B:=V^*-D$$
where $D$ consists of $\sigma:=(\sigma_1,..,\sigma_s)\in V^*$  such that $Y_\sigma=Y_{\sigma_1,..,\sigma_s}$ is not smooth of codimension $s$. Let $\cY$ denote the family of complete intersections $Y_\sigma$ parameterized by $B$.
\item{\bf X3.} By {\bf X1}, we have canonical $G$-equivariant maps $X\ra\PP V_i$. Denote the map by
$$
\phi: X\ra\PP V_1\times\cdots\times\PP V_s
$$
and let $\hat X$ be the cone over the image in $V=V_1\times\cdots\times V_s$. Since the $L_i$ are $G$-equivariant,
we have a holomorphic representation $G\ra\Aut V$. We extend this to a representation
$$
\hat G=G\times(\CC^\times)^s\ra\Aut V
$$
by letting the $i$th copy of $\CC^\times$ in $(\CC^\times)^s$ act on $V_i$ by the usual scaling.
\end{itemize}


The next proposition gives the most important case in which {\bf X1-X3} hold.

\begin{prop}
Suppose $L=L_1=-K_X$ ($s=1$) is very ample. Then {\bf X1-X3} hold for the family $\cY$ of CY hypersurfaces in $X$ with $H=\C^\times$.
\end{prop}
\begin{proof}
Obviously $H^0(L+K_X)\neq 0$, and $H^0(L)$ is base point free. By Theorem \ref{Calabi-form}, $\C^\times-K_X^\times\ra X$ is a CY bundle. By Theorem \ref{KX-theorem}, $L=-K_X=L_{\chi_1^{-1}}$. So, {\bf X1} holds. Since $L$ is very ample, {\bf X2} also follows. 
\end{proof}

\begin{rem}
Since each $L_i$ in {\bf X1} is nontrivial, and since $H^0(L_i)$ is base point free, it has dimension at least 2. In a typical application in this paper, the $L_i$ are very ample, in which case the $H^0(L_i)$ are base point free and {\bf X2} hold automatically. It should also be noted that if $H-M\ra X$ is a principal bundle such that the natural map
$$
\Hom(H,\C^\times)\ra \Pic_G(X)
$$
is onto (cf, Proposition \ref{Popov}), then the assumption that $L=L_\chi$ in {\bf X1} is redundant. 
\end{rem}

Let $\tau\in H^0(L+K_X)$, and $\sigma\in B$. Then by Theorems  \ref{residue} and \ref{Calabi-form} (in the case $M=K_X^\times$, $H=\CC^\times$, $\Omega=\iota_{x_0}\hat\omega$)
$$
\tau\Omega_\sigma:={\tau\iota_{x_0}\hat\omega\over \sigma_1\cdots \sigma_s} 
$$
defines a meromorphic $d$-form on $X$ with pole along the divisor $\cup_{i=1}^s Y_{\sigma_i}$.

\begin{thm}\label{globalPR} (Global Poincar\'e residue) Assume {\bf X1-X3} with $H=\C^\times$. Then the global Poincar\'e residue map ${\bf R}:H^0(L+K_X)\ra H^0(B,\HH^{top})$ is $G$-equivariant and is given by
\begin{equation}\label{GPR}
\R_\sigma(\tau)=\Res\tau\Omega_\sigma.
\end{equation}
\end{thm}
\begin{proof}
We verify eqn. \ref{GPR} first. By definition, for fixed $\sigma\in B$, the restriction map $\R_\sigma$ is induced on global sections by the sheaf homomorphism \cite{GH}
$$
\cO(L+K_X){\br\Psi_\psi\over\cong} \Omega_X^d(\cup_i Y_{\sigma_i})~{\br PR\over\ra}~\Omega^{d-s}_{Y_\sigma}=\cO(K_{Y_\sigma})
$$
where $\Psi_\psi:\cO(L+K_X)=\cO(L)\otimes\cO(K_X)\ra\Omega_X^d(\cup_i Y_{\sigma_i})$ is given by $\phi\otimes\omega\mapsto {\phi\omega\over\psi}$ with $\psi=\sigma_1\cdots\sigma_s$, and $PR$ is the (codimension $s$) Poincar\'e residue map. So,
$$
\R_\sigma=\Res\circ\Psi_\psi:H^0(L+K_X)\ra H^0(K_{Y_{\sigma}}).
$$
For eqn. \ref{GPR}, it is enough to show that for any given nonzero section $\psi\in H^0(L)$, $\Psi_\psi$ on global sections is given by
$$
H^0(L+K_X)\ra H^0(\Omega_X^d(Y)),~~\tau\mapsto{\tau\iota_{x_0}\hat\omega\over\psi}
$$
where $Y$ is the divisor $\psi=0$ in $X$. This map is well-defined by Theorems \ref{residue} and \ref{Calabi-form} (in the case $M=K_X^\times$, $H=\CC^\times$, $\Omega=\iota_{x_0}\hat\omega$.) 

Since $L=L_\chi$ and $\chi$ lies in the cyclic group $\Hom(H,\CC^\times)$ generated by $\chi_1$, it follows that $L=l K_X$ for some $l\in\ZZ$. Fix a coordinate $(U,w)$ on $X$ and put $dw:=dw_1\wedge\cdots\wedge dw_d$. Then $\tau\in H^0(L+K_X)$ and $\psi\in H^0(L)$ have the local form on $U$
$$
\tau=\tau^U(w) dw^{\otimes(l+1)},\blank\psi=\psi^U(w) dw^{\otimes l}.
$$
Since $K_X=M\times_H\C_{\chi_1}$ is trivialized on $U$ by 
$$
M|U\times_H\C_{\chi_1}\ra U\times\C,~~~[\eta,1]\mapsto([\eta],z_w(\eta))
$$
it follows from Proposition \ref{compare-sections} that as functions on $M|U$, $\tau,\psi$ are given by 
$$
\tau|U=\tau^U(w) z_w^{-l-1},\blank\psi|U=\psi^U(w) z_w^{-l}.
$$
Thus by Theorem \ref{Calabi-form}, we get
$$
{\tau\iota_{x_0}\hat\omega\over\psi}|U={\tau^U dw\over\psi^U}.
$$
On the other hand, for $\phi\in\cO(L|U)$ and $\omega\in\cO(K_X|U)$, we have
$$
\Psi_\psi:\phi\otimes\omega\mapsto{\phi^U\omega^U dw\over\psi^U}
$$
Since $\tau|U$ is a finite sum of tensors of the form $\phi\otimes\omega$, it follows that $\Psi_\psi(\tau|U)$ coincides with ${\tau\iota_{x_0}\hat\omega\over\psi}|U$ above. This shows that
$$
\Psi_\psi(\tau)={\tau\iota_{x_0}\hat\omega\over\psi}
$$
proving eqn. \ref{GPR}.

Let $g\in\Aut X$. Then $g$ induces isomorphisms
\begin{eqnarray*}
& &g_\cup:H^0(\Omega_X^d(\cup_i Y_{\sigma_i}))\ra  H^0(\Omega_X^d(\cup_i Y_{\Gamma_g\sigma_i}))\cr
& &g_\cap:H^0(\Omega_{Y_\sigma}^{d-s})\ra  H^0(\Omega^{d-s}_{Y_{\Gamma_g\sigma}})\cr
\end{eqnarray*}
such that $g_\cap\circ\Res =\Res\circ g_\cup$.
\comment{See my NTU hand written lecture notes for details.}
So, for $\tau\in H^0(L+K_X)$, we have
\begin{eqnarray*}
\Gamma_g\R(\tau)&=&g_\cap \Res\tau\Omega_{\Gamma_g^{-1}\sigma}
=\Res g_\cup (\tau\Omega_{\Gamma_g^{-1}\sigma})\cr
&=&\Res (\Gamma_g\tau)\Omega_\sigma=\R(\Gamma_g\tau).
\end{eqnarray*}
This shows that ${\bf R}$ is equivariant.
\end{proof}



\section{Comparing global Poincar\'e residue formulas}

We now give two more applications of the canonical form. We have seen that the CY bundle $\CC^\times -K_X^\times\ra X$ gives us a formula for the global Poincar\'e residue $\R$ (Theorem \ref{globalPR}.) However, the formula is subject to the condition that the bundle $L$ in assumption {\bf X1} has the form $L_\chi$, where $\chi$ is an $\C^\times$-character. By Lemma \ref{L-star}, this restricts $L$ to a power of $K_X$ (cf. Example \ref{Kx-Ominus1}.) We now show that $\R$ can be realized by using any other CY bundle $H-M\ra X$. This introduces enormous flexibility to the use of the global Poincar\'e residue. For example, when $X$ is a homogeneous manifold $X=G/P$, we will see that the natural bundle
$$
P-G\ra X
$$
has a CY structure. This allows us to realize {\it all} line bundles on $X$ in terms of $P$-characters, by Proposition \ref{Popov}. So, for any homogeneous space, the assumptions in {\bf X1} is superfluous. The same is true for any toric manifold.  The next example, which lives in both worlds (i.e. homogeneous and toric), shows two CY $\C^\times$-bundles that realize $\R$ in two different ways.

\begin{exm}\label{Kx-Ominus1}
Projective hypersurfaces.
\end{exm}
\cut

Let $X=\PP^d$. According to Lemma \ref{L-star}, the CY bundle $\C^\times-K_X^\times\ra X$ allows us to realize the global Poincar\'e residue by means of Theorem \ref{globalPR}, for the line bundle $L=\cO(k)$ ($k>0$) only when $L=lK_X=\cO(-ld-l)$ for some $l\in\ZZ$. This forces $k$ to be an integer multiple of $(d+1)$. By contrast, the bundle $\C^\times-\cO(-1)^\times\ra X$, which is a CY bundle by Theorem \ref{obstruction}, can realize the global Poincar\'e residue for {\it all} $\cO(k)$, because each $\cO(k)$ is a power of $\cO(-1)$ (cf. Lemma \ref{L-star}.) Note that by viewing $X$ as a homogeneous space of $SL_{d+1}$ gives a third global Poincar\'e residue formula, for all line bundle $L$ as well.

We now return to the context of \ts\ref{sec-globalPR} beginning with the assumptions {\bf X1-X3}. We will prove that the global Poincar\'e residue formula derived from any given CY bundle, by means of Theorem \ref{residue}, agrees with $\R$ of Definition \ref{period-integrals}. Given a CY bundle $H-M\ra X$, an $H$-character $\chi$, and a section $\psi\in H^0(L_\chi)$, let $\psi_M:M\ra\C$ denote the $\chi$-equivariant function representing $\psi$ (cf. Proposition \ref{function-section}.)


\begin{thm} (Uniqueness of global Poincar\'e residue)\label{unique-globalPR}
Let $H-M\ra X$ be a bundle with the CY structure $(\omega_M,\chi_M)$. Assume {\bf X1-X3}.
Fix independent vector fields $x_i$ generated by $H$ on $M$. There is a constant $c$ such that, for $\tau\in H^0(L_\chi+K_X)$, $\sigma\in B$, and $\psi:=\sigma_1\cdots\sigma_s\in H^0(L_\chi)$, we have
$$
\R_\sigma(\tau)=c\Res{\tau_M\over\psi_M}\iota_{x_1}\cdots\iota_{x_q}\omega_M.
$$
Moreover, if $\cup_{\sigma\in B} Y_\sigma$ is dense in $X$, then $\R$ is injective.
\end{thm}

\comment{On the last assertion, see subsection Injectivity \& surjectivity of R, dual variety, discriminant locus of flag.tex}

For the first assertion, the overall strategy is about the same as proof of Theorem \ref{globalPR}, but there are some important differences. One is that the local frames for line bundles used in the former proof are too restrictive to handle a general CY bundle $M$. Another difference is that in present setting, we will need to invoke the uniqueness theorem for CY structures on $M$. We were able to avoid this earlier essentially because the canonical form was particularly simple in local coordinates, which in turn made it a lot easier to identify the CY structure in question. We now begin with some preparations.

Just as in Theorem \ref{globalPR}, since $\R_\sigma=\Res\circ\Psi_\psi$, it suffices to show that we can choose a normalization for $\omega_M$ so that
\begin{equation}\label{to-prove}
\Psi_\psi(\tau)={\tau_M\over\psi_M}\iota_{x_1}\cdots\iota_{x_q}\omega_M
\end{equation}
as meromorphic forms on $M$. Note that the left side is a meromorphic form on $X$, but can be regarded equivalently as an $H$-basic form on $M$.

We shall prove that eqn. \ref{to-prove} holds on $M|U$ for all sufficiently small open set $U\subset X$. Fix a local trivialization of the $H$-bundle $M$, i.e. an $H$-equivariant map of the form
$$
M|U\ra U\times H,~~m\mapsto([m],\alpha(m))
$$
which we will label $\alpha$. We will also use the letter $\beta$ to denote such local trivializations of the bundle $M$. That the map is $H$-equivariant implies that
$$
\alpha(mh^{-1})=\alpha(m)h^{-1}\blank(h\in H.)
$$
It follows that for any given $H$-character $\xi$, we have
\begin{equation}\label{e3}
\xi(\alpha(mh^{-1}))=\xi(\alpha(m))\xi(h)^{-1}.
\end{equation}
By Proposition \ref{compare-sections}, this defines a unique local trivialization of $L_\xi$
$$
L_\xi|U\ra U\times\C,~~~[m,1]\mapsto([m],\xi(\alpha(m)))
$$
which we will also label $\alpha$. By the same proposition, for each section $\psi\in H^0(L_\xi)$, we have
\begin{equation}\label{e4}
\psi_M(m)=\psi^\alpha([m])\xi(\alpha(m))^{-1}\blank(m\in M|U.)
\end{equation}
The main point here is that a single local trivialization $\alpha$ of the $H$-bundle $M$ determines simultaneously local trivializations for {\it all} line bundles $L_\xi$ in a canonical way. In turn, such a local trivialization determines a local frame $e_\alpha$ for $L_\xi|U$, and relates via eqn. \ref{e4}, the $\xi$-equivariant function $\psi_M\in\cO(M)_\xi$ representing a section of $L_\xi$, and its local expression on $U\subset X$:
$$
\psi|U=\psi^\alpha e_\alpha.
$$

For the proof of eqn. \ref{to-prove}, we apply this to the line bundle $L=L_\chi$ named in the theorem, and the canonical bundle  $K_X=L_{\chi_K}$, with $\chi_K:=\chi_M\chi_\gh$. Thus $\alpha$ determines their local frames, which we denote by $e_\alpha$ and $\varepsilon_\alpha$ respectively. In particular, the sections $\tau\in H^0(L+K_X)$, $\psi\in H^0(L)$ have the local expressions
$$
\tau|U=\tau^\alpha e_\alpha\otimes\varepsilon_\alpha,\blank\psi|U=\psi^\alpha\varepsilon_\alpha.
$$
Likewise we have similar expressions corresponding to any other local trivialization $\beta$ of $M$.

Next, we unwind both sides of eqn. \ref{to-prove} to reduce it further. Recall that the sheaf map
$$
\Psi_\psi:\cO(L)\otimes\cO(K_X)\ra\Omega^d_X(Y)
$$
(where $Y$ is the divisor $\psi=0$) is given by
$$
a\otimes b=a^\alpha e_\alpha\otimes b^\alpha\varepsilon_\alpha\mapsto {a^\alpha b^\alpha\over\psi^\alpha}\varepsilon_\alpha.
$$
Since $\tau^\alpha$ is locally a finite sum of terms of the form $a^\alpha b^\alpha$, it follows that on $U$,
$$
\Psi_\psi(\tau)={\tau^\alpha\over\psi^\alpha}\varepsilon_\alpha.
$$

On the other hand, by eqns. \ref{e3}-\ref{e4}, we have on $M|U$,
$$
{\tau_M\over\psi_M}
={\tau^\alpha\over\psi^\alpha}{(\chi\circ\alpha)^{-1}(\chi_K\circ\alpha)^{-1}\over(\chi\circ\alpha)^{-1}}
={\tau^\alpha\over\psi^\alpha}(\chi_K\circ\alpha)^{-1}
$$
So, we have reduce proving eqn. \ref{to-prove} to
\begin{equation}\label{e5}
(\chi_K\circ\alpha)\varepsilon_\alpha=\Omega:=\iota_{x_1}\cdots\iota_{x_q}\omega_M.
\end{equation}

\begin{lem} 
$(\chi_K\circ\alpha)\varepsilon_\alpha$, as a $d$-form on $M|U$, is independent of the choice of local trivialization $\alpha$ of $M$ over $U$. Thus together the $\{(\chi_K\circ\alpha)\varepsilon_\alpha\}$ define a global nowhere vanishing $d$-form $v$ on $M$.
\end{lem}
\begin{proof}
Again, that $(\chi_K\circ\alpha)\varepsilon_\alpha$ is nowhere vanishing is obvious, since $\chi_K$ takes value in $\C^\times$ and $\varepsilon_\alpha$ is a frame of $K_X$ (but viewed as an $H$-basic form on $M|U$.)

By Proposition \ref{compare-sections}, for any local trivializations $\alpha,\beta$ of $M$ over $U$, we have
\begin{eqnarray*}
\tau^\alpha(\chi\circ\alpha)^{-1}(\chi_K\circ\alpha)^{-1}&=&\tau^\beta(\chi\circ\beta)^{-1}(\chi_K\circ\beta)^{-1}\cr
\psi^\alpha(\chi\circ\alpha)^{-1}&=&\psi^\beta(\chi\circ\beta)^{-1}\cr
\end{eqnarray*}
It follows that
$$
{\chi_K\circ\alpha\over\chi_K\circ\beta}=
{\varphi^\alpha\over\varphi^\beta}
$$
where $\varphi^\alpha:={\tau^\alpha\over\psi^\alpha}$. Note also that on $M|U$, we have
$$
\varphi^\alpha\varepsilon_\alpha=\Psi_\psi(\tau)=\varphi^\beta\varepsilon_\beta.
$$
It follows that
$$
(\chi_K\circ\alpha)\varepsilon_\alpha
=(\chi_K\circ\beta){\varphi^\alpha\over\varphi^\beta}\varepsilon_\alpha
=(\chi_K\circ\beta){\varphi^\beta\over\varphi^\beta}\varepsilon_\beta
=(\chi_K\circ\beta)\varepsilon_\beta.
$$
This completes the proof.
\end{proof}

Recall that $x_1,..,x_q$ are independent global vector fields generated by the $H$-action on $M$. Fix a local frame $y_1,...,y_d$ of the bundle $TU\subset T(U\times H)$ once and for all (independent of $\alpha$). Since we have an $H$-equivariant isomorphism
$$
M|U\ra U\times H,~~~m\mapsto([m],\alpha(m))
$$
and since $H$ acts only on the $H$ factor on $U\times H$, the vector fields
$\alpha_*x_i$ are clearly independent of the $y_k$ on $U\times H$. It follows that there exists unique 1-forms $\omega_j^\alpha$, $j=1,...,q$, on $U\times H$ such that
$$
\iota_{\alpha_*x_i}\omega_j^\alpha=\delta_{ij},\blank\iota_{y_k}\omega^\alpha_j=0.
$$
These equations imply that the $\omega_j^\alpha$ are local sections of $T^*H\subset T^*(U\times H)$. Put
$$
y_k^\alpha=(\alpha_*)^{-1}y_k,\blank\omega_j^\alpha\equiv\alpha^*\omega_j^\alpha.
$$
Then on $M|U$, we have
$$
\iota_{x_i}\omega_j^\alpha=\delta_{ij},\blank\iota_{y_k^\alpha}\omega^\alpha_j=0.
$$
Put
$$
u^\alpha:=\omega_q^\alpha\wedge\cdots\wedge\omega_1^\alpha
$$
as a $q$-form on $M|U$

\begin{lem}
The form $u^\alpha\wedge v$ is independent of the choice of local trivialization $\alpha$ of $M$ over $U$. Thus together the $\{u^\alpha\wedge v\}$ define a global nonvanishing $(q+d)$-form $\omega'$ on $M$.
\end{lem}
\begin{proof}
Since $v$ is $H$-horizontal, it follows that $\iota_{x_1}\cdots\iota_{x_q}(u^\alpha\wedge v)=v$, which is nonvanishing. It follows that
 $u^\alpha\wedge v$ is nonvanishing.

Suppose $\beta$ is another local trivialization of $M$ over $U$. Then there exists holomorphic functions $\lambda_{\alpha\beta}$ such that
$$
\blank u^\alpha\wedge v =\lambda_{\alpha\beta} u^\beta\wedge v.
$$
Applying $\iota_{x_1}\cdots\iota_{x_q}$ (which is globally defined on $M$, i.e. independent of $\alpha$ and $U$) to both sides, we find that $\lambda_{\alpha\beta}=1$.
\end{proof}

\begin{lem}
The nonvanishing form $\omega'\in H^0(K_M)$ in the preceding lemma is 
a CY structure of $M$.
\end{lem}
\begin{proof}
 Let $h\in H$. We have on $U$,
$$
\Gamma_h u^\alpha=\Gamma_h(\omega_q^\alpha\wedge\cdots\wedge\omega_1^\alpha)=\chi_\gh(h)^{-1}u^\alpha
$$
since $u^\alpha$ is a trivializing section of the dual bundle of $M|U\times_H\wedge^q\gh$ (and the latter transforms by $\chi_\fh(h)$.)
We also have (keeping in mind that $h$ acts trivially on $T^*U$ of which $\varepsilon_\alpha$ is a section)
$$
(\Gamma_h v)(m)=v(mh)=\chi_K(\alpha(mh))\varepsilon_\alpha([mh])=\chi_K(h)v(m).
$$
So, $\Gamma_h v=\chi_K(h)v$. Since $\chi_K=\chi_M\chi_\gh$, it follows that
$$
\Gamma_h(u^\alpha\wedge v)=\chi_M(h)(u^\alpha\wedge v).
$$
This completes the proof.
\end{proof}

\begin{proof}
We are now ready to prove eqn. \ref{e5}, and to complete the proof of Theorem \ref{unique-globalPR}. By Theorem \ref{uniqueness} and the preceding lemma, it follows that $\omega'=u^\alpha\wedge v$ is equal $\omega_M$, after a suitable scalar normalization.
Finally, this implies that
$$
v=\iota_{x_1}\cdots\iota_{x_q}(u^\alpha\wedge v)=\iota_{x_1}\cdots\iota_{x_q}\omega_M=\Omega
$$
which is the desired eqn. \ref{e5}. 

We now prove the last assertion of Theorem \ref{unique-globalPR}. Suppose that $\R(\tau)=0$. This is equivalent to that $Y_\sigma\subset Z_\tau$ for all $\sigma\in B$, where $Z_\tau$ is the zero locus of $\tau$. Since $\cup_{\sigma\in B} Y_\sigma$ is dense in $X$, it follows that $Z_\tau$ must be all of $X$. Hence $\tau=0$.
\end{proof}

\begin{rem}
On the last assertion, the assumption that $\cup_{\sigma\in B}Y_\sigma$ is dense in $X$ is superfluous in the codimension $s=1$ case. Hence in this case, $\R$ is always in injective. This is also true for $s\geq1$ if one assumes that the line bundles $L_1,..,L_s$ are very 
ample.
\end{rem}

\begin{defn}
Given a CY bundle $H-M\ra X$, we shall say that a CY structure $\omega_M$ on $M$ is normalized if the constant $c$ in the preceding theorem is 1. In other words
$$
\R_\sigma(\tau)=\Res \tau\Omega_\sigma.
$$
for all $\sigma\in B$ and $\tau\in H^0(L+K_X)$. Note that this notion is independent of the choice of $L$.
\end{defn}

\begin{exm}
Projective hypersurfaces, revisited.
\end{exm}
\cut

Let $X=\PP^d$, and $L=\cO(n)$ with $n\geq d+1$. Then $L+K_X=\cO(n-d-1)$ and $H^0(L+K_X){\br\Psi_\sigma\over\cong }H^0(\Omega_X^d(Y_\sigma))$ for any general hypersurface $Y_\sigma$ of degree $n$ in $X$. 
As shown in Example \ref{Kx-Ominus1}, we can use the CY bundle $\C^\times - \cO(-1)^\times\ra X$ to realize $L$ as $L_\chi$, for some $\C^\times$-character $\chi$. By Proposition \ref{PR-isomorphism}, $\dim H^0(L+K_X)=h^{d-1,0}(Y_\sigma)$ which is the rank of $\HH^{top}$. 
By the preceding theorem, the sections in the image of $\R$ provide a global trivialization of the vector bundle $\HH^{top}$ in this case.

\section{Tautological systems for period integrals: general case}\label{sec-tautological}


Let $\hat G$ be a complex Lie group, and $\rho:\hat G\ra\Aut V$ be a finite dimensional holomorphic representation such that $\CC^\times\subset\rho(\hat G)$. Let 
$$
\hat\fg\ra\End V, ~~x\mapsto Z_x
$$
be the corresponding Lie algebra homomorphism. Since $\CC[V^*]=\Sym V$, we will view $\End V$ as a Lie subalgebra of $\Der\CC[V^*]$. Thus $Z_x\in\Der \CC[V^*]$ for $x\in\hat\fg$. Let $\hat X\subset V$ be a given $\hat G$-invariant subvariety in $V$, and $I(\hat X,V)\subset\CC[V]$ be the vanishing ideal of $\hat X$. Note that since $\hat X$ is invariant under $\CC^\times\subset\rho(\hat G)$, the ideal is homogeneous.

Since we have a canonical symplectic form $\bra,\ket$ on $T^*V=V\times V^*$, each linear function $\zeta\in V^*$ uniquely defines a derivation $\partial_\zeta\in \Der\CC[V^*]$, by the formula 
$$
\partial_\zeta a=\bra a,\zeta\ket,~~ a\in V.
$$
Let $D_{V^*}$ be the algebra of global differential operators on $V^*$.

\begin{defn} Fix a Lie algebra homomorphism $\beta:\hat\fg\ra\C$. The tautological system 
$$\tau(\hat X,V,\hat G,\beta)$$ is the $D_{V^*}$-module $D_{V^*}/J$ where $J$ is the left ideal of $D_{V^*}$ generated by the following operators: $\{p(\partial_\zeta)|p(\zeta)\in I(\hat X,V)\}$, $\{Z_x+\beta(x)| x\in\hat\fg\}$. 
\end{defn}

We call generators of the form $Z_x+\beta(x)$, the {\it symmetry operators or $\hat G$-operators}, and those of the form $p(\partial_\zeta)$, the {\it embedding polynomial operators or polynomial operators}.
We note that the notion of a tautological system here is slightly broader than that given in \cite{LSY1}, in that the formulation here is purely algebraic and does not require that one starts with a $G$-manifold, or that $G$ is reductive, or that $\beta$ comes from a character of $\C^\times$.



Using dual bases, $a_i$ and $\zeta_i$ of $V,V^*$, we can write an element of $ \CC[V]=\CC[\zeta_1,\zeta_2,...]$ as a polynomial $p(\zeta)=p(\zeta_1,\zeta_2,...)$, and $p(\partial_\zeta)$ as a partial differential operator $p({\partial\over\partial a_1},{\partial\over\partial a_2},...)$ with constant coefficients, acting on functions of the variables $a_1,a_2,...$. If $(x_{ji})$ is the matrix representing $x\in\hat\fg$ acting on $V$ in the basis $a_i$, i.e. $x\cdot a_i=\sum_j x_{ji}a_j$, then 
$$Z_x=\sum x_{ji}a_j{\partial\over\partial a_i}.$$

Let $j:V\into W$ be a $G$-module homomorphism, $\pi:W^*\onto V^*$ the dual map. This induces an $G$-equivariant map on structure sheaves $\pi^\#:\cO_{V^*}\ra\pi_*\cO_{W^*}$, $f\mapsto f\circ \pi$ for $f\in\cO_{V^*}(U)$, and the induced homomorphism on germs is also a $\fg$-module homomorphism:
\begin{equation}\label{G-operator}
(Z_x^Vf)\circ\pi=Z^W_x(f\circ\pi)
\end{equation}
where $f\in\cO_{V^*}$.  Here $Z_x^V,Z_x^W$ are the $G$ operators on $V^*,W^*$ respectively. Likewise for the Euler operators (with the same character $\beta$.) 

Now $\pi:W^*\onto V^*$ induces the $G$-equivariant algebra homomorphisms $\pi:\CC[W]\onto\CC[V]$ and
$$
\pi:\C[\partial_{\zeta_W}|\zeta_W\in W^*]\onto
\C[\partial_{\zeta_V}|\zeta_V\in V^*]~~~\partial_{\zeta^W}\mapsto\partial_{\pi\zeta^W}.
$$
It is straightforward to check that for $f\in\cO_{V^*}$, $p(\partial_{\zeta^W})\in \C[\partial_{\zeta^W}]$, we have
\begin{equation}\label{poly-operator}
[(\pi p(\partial_{\zeta^W})) f]\circ\pi=p(\partial_{\zeta^W})(f\circ\pi).
\end{equation}

Let $\cS_{V^*}$ (likewise $\cS_{W^*}$) be the subsheaf of $\cO_{V^*}$ whose stalks consists of germs annihilated by the defining ideal $J$ of the D-module $\cM=\tau(X,\phi_V,G,\beta)$. Then we have a canonical isomorphism (of $\CC_{V^*}$-modules) from $\cS_{V^*}$ to the solution sheaf $\cH om_{\cD_{V^*}}(\cM,\cO_{V^*})$ of $\cM$. Under this identification, we can therefore view $\cS_{V^*}$ as the solution sheaf of $\cM$.



\begin{lem} (Change of variables)\label{change-of-variables}
Let $j:V\into W$ be a $\hat G$-module homomorphism and $\pi:W^*\onto V^*$ its dual. Let $\cS_{V^*}\subset\cO_{V^*}$ and $\cS_{W^*}\subset\cO_{W^*}$ be the solution sheaves of the D-modules $\tau(\hat X,V,\hat G,\beta)$ and $\tau(\hat X,W,\hat G,\beta)$. Then $\pi^\#$ maps $\cS_{V^*}$ isomorphically onto $\cS_{W^*}$.
\end{lem}

The proof in \cite{LSY1} carries over with little change.
The lemma allows us to give different descriptions to essentially the same D-module, by choosing different $\hat G$-modules as targets. The net effect of changing the target from $V$ to $W$ in the initial data of our tautological system, is that we introduce additional linear variables, and at the same time, additional first order operators corresponding to the linear forms in $V^\perp\subset W^*$. 

\begin{thm}  \cite{LSY1}\label{holonomic}
Assume that the $\hat G$-variety $\hat X$ has only a finite number of $\hat G$-orbits. Then the tautological system $\tau(\hat X,V,\hat G,\beta)$ is holonomic. In particular, the space of solutions at a generic point is finite dimensional.
\end{thm}


\begin{exm} Homogeneous $G$-variety.
\end{exm}
\cut

Let $G$ be a semisimple group, and $\beta:\CC\ra\CC$ a linear function.
By the Borel-Weil theorem, for any ample line bundle $L$ on $X=G/P$, $V=H^0(X,L)^*$ is an irreducible $G$ module, say, of highest weight $\lambda$. It defines a $G$-equivariant embedding $X\into\PP V$. Let $\hat X$ be the cone over the image of $X$, and let $\CC^\times$ acts on $V$ by the usual scaling. Put $\hat G=G\times\CC^\times$, and extend $\beta$ to $\hat\fg=\fg\oplus\C\ra\C$ by setting $\beta(\fg)=0$.  By a theorem of Kostant and Lichtenstein, the ideal $I(\hat X,V)$ is generated by quadratic forms. In particular, aside from the symmetry operators  $Z_x+\beta(x)$ ($x\in\hat\fg$), all other generators of our tautological system $\tau(\hat X,V,\hat G,\beta)$ are second order. The latter can be more explicitly described as follows. Let
$$
C=\sum_{\alpha>0}(x_\alpha x_{-\alpha}+x_{-\alpha}x_\alpha)+\sum_{i=1}^{rk(G)} h_i^2
$$
be the second Casimir operator in the enveloping algebra $U\fg$. We have the following theorem:

\begin{thm}\label{Lichtenstein}\cite{Kostant} \cite{Lich} $I(\hat X,V)$ is generated by $\half m_\lambda(m_\lambda+1)$ ($m_\lambda=\dim V^*$) quadrics, given by
$$C(u\otimes u)-\bra 2\lambda+2\delta,2\lambda\ket(u\otimes u)\blank(u\in V^*)$$
where $\delta$ is half the sum of the positive roots of $\fg$.
\end{thm}

Note that the quadrics above are viewed as elements of $Sym^2 V^*\subset\CC[V]$. Later, we will give another description of this D-module, by embedding $V$ into a natural but larger $\hat G$-module $W$, using the classical Veronese and Segre (and the incidence) maps. 

\begin{exm} Schubert varieties.
\end{exm}
\cut

Let $G$ be a semisimple group, $P$ a parabolic and $B\subset P$ a Borel subgroup of $G$.
A Schubert variety $X$ is the closure of a $B$-orbit in $G/P$. Since $X$ has only finite number of $B$-orbits, Theorem \ref{holonomic} applies to the tautological system $\tau(\hat X,W,\hat B,\beta)$ ($\hat B=B\times\CC^\times$) for any $B$-equivariant embedding $\phi:X\into\PP W$. Here $\hat X$ is the cone over $\phi(X)$ in $W$. This class of D-modules will be investigated in a future paper.

\begin{exm} Invariant theory.
\end{exm}
\cut

Next, we point out that the problem of describing the solution sheaf of a tautological system can be thought of as a generalization of a problem in classical invariant theory. Let $G$ be a connected semisimple algebraic group and $V$, a $G$-module. Consider the D-modules $\tau(V,V,G\times\CC^\times,\beta)$, where $\beta:\fg\oplus \CC\ra\CC$ with $\beta(\fg)=0$ and $\beta(1)\in\ZZ_<$. Since the ideal $I(V, V)$ is trivial, the only generators of our tautological system are the Euler operator $\sum_i a_i{\partial\over\partial a_i}+\beta(1)$, and the symmetry operators corresponding to $G$. Thus a function $f\in\cO_{V^*}(U)$ is a solution of the D-module iff it is $G$-invariant and is homogeneous of degree $-\beta(1)$. If one demands that $f$ be a global polynomial solution on $V^*$, then $f$ is nothing but an element of $\CC[V^*]^G$ of a given degree.

To construct differential equations for period integrals, Definition \ref{period-integrals}, we now return to the context of \ts\ref{sec-globalPR} and assumptions {\bf X1-X3}. Thus we fix a CY bundle
$$
H-M\ra X
$$
with a normalized CY structure. First, by Theorem \ref{unique-globalPR}, the period integrals of the family $\cY$ of complete intersections in $X$ are given by
$$
\Pi_\gamma(\tau;\sigma):=\int_\gamma\R_\sigma(\tau)=\int_\gamma \Res \tau_M\Omega_\sigma,\blank\gamma\in H_{d-s}(Y_\bullet,\ZZ)
$$
viewed as locally defined functions in the variable $\sigma\in B$, with $\tau\in H^0(L+K_X)$ fixed.

\begin{thm}\label{general-tautological} 
Let $H-M\ra X$ be a CY bundle and assume {\bf X1-X3}. Let $\gamma\in H_{d-s}(Y_\bullet,\ZZ)$, and $\tau_0\in H^0(L+K_X)$ be a fixed eigenvector of $\hat\fg$: $x\cdot\tau_0=\lambda_0(x)\tau_0$ for $x\in\hat\fg$ for some $\lambda_0\in{\hat\fg}^*$.
Then the period integrals $\Pi_\gamma(\tau_0;\sigma)$ of the family $\cY$ are solutions to the tautological system $\tau(\hat X,V,\hat G,\beta)$, where $\beta\in{\hat\fg}^*=\fg^*\oplus\CC^s$ is the vector $\beta=(\lambda_0;1,1,..,1)$.
\end{thm}
\begin{proof}
Let $p(\zeta)\in I(\hat X,V)$. Since $\hat X$ is $\hat G$-invariant subvariety of 
$$
V=V_1\times\cdots\times V_s
$$ 
it is invariant under scaling by $(\CC^\times)^s$. In particular, we may as well assume that $p(\zeta)$ is homogeneous with respect to this scaling, say of degrees $(k_1,...,k_s)$. Now if $\zeta_1,..,\zeta_m$ form a basis of $V_i^*$ and $a_1,..,a_m$ form the dual basis, then for $\zeta\in V_i^*$, $i$ fixed, we have $\partial_\zeta\sum_j a_j\zeta_j=\sum_j\bra a_j,\zeta\ket\zeta_j=\zeta$, hence
$$\partial_\zeta {1\over(\sum a_j\zeta_j)^l}=-l\zeta {1\over (\sum a_j\zeta_j)^{l+1}}.$$
It follows that
$$
p(\partial_\zeta)\tau_0\Omega_\sigma=(-1)^{k_1+\cdots+k_s} k_1!\cdots k_s! p(\zeta){\tau_0\Omega\over \sigma_1^{k_1+1}\cdots\sigma_s^{k_s+1}}
$$
where
$$\Omega=\iota_{x_1}\cdots\iota_{x_q}\omega_M
$$
and $\omega_M$ is the normalized CY structure on $M$. Since $p(\zeta)=0$ on $\hat X$, it follows that $p(\partial_\zeta)$ annihilates $\Pi_\gamma(\tau_0;\sigma)$.

We can write
$$
\Pi_\gamma(\tau_0;\sigma)=\int_{\tau(\gamma)}\tau_0\Omega_\sigma
$$
where $\tau(\gamma)$ is a tube over the cycle $\gamma$. Let $g\in\Aut X$. Since the period is the Poincar\'e pairing $\bra\tau(\gamma),\tau_0\Omega_\sigma\ket$ on $X-\cup_i Y_{\sigma_i}$, it is invariant under $g$:
$$
\bra\tau(\gamma),\tau_0\Omega_\sigma\ket=\bra (g_*)^{-1}\tau(\gamma),g^*(\tau_0\Omega_\sigma)\ket.
$$
Now let $g\in G$ be close to identity. Then $(g_*)^{-1}\tau(\gamma)=\tau(\gamma)$. By Theorem \ref{residue}, $\Omega$ is $G$-invariant, so 
$$
\bra \tau(\gamma),{\tau_0\Omega\over\sigma_1\cdots\sigma_s}\ket=\bra\tau(\gamma),g^*({\tau_0\over\sigma_1\cdots\sigma_s})\Omega\ket.
$$
For $x\in\fg$, consider the action of the 1-parameter subgroup $g=g_t=exp(tx)$ of $G$. We have
\begin{eqnarray*}
{d\over dt}|_{t=0} g_t^*({\tau_0\over\sigma_1\cdots\sigma_s})&=&{x\cdot\tau_0\over\sigma_1\cdots\sigma_s}-\sum_i{\tau_0~x\cdot \sigma_i\over \sigma_1\cdots\sigma_i^2\cdots\sigma_s}\cr
&=&{\lambda_0(x)\tau_0\over\sigma_1\cdots\sigma_s}+Z_x({\tau_0\over\sigma_1\cdots\sigma_s})\cr
\end{eqnarray*}
 It follows that
$$
0=\bra\tau(\gamma),(Z_x+\lambda_0(x))\tau_0\Omega_\sigma\ket=(Z_x+\lambda_0(x)) \bra \tau(\gamma),\tau_0\Omega_\sigma\ket
$$
Finally, fix $i$ and let $y=(0;0,..,1,..,0)\in\fg\oplus\CC^s=\hat\fg$, where the 1 is in the $i$th slot of $y$. Then $Z_y$ is the Euler operator $\sum_{j=1}^m a_j{\partial\over\partial a_j}$. Since ${1\over \sigma_i}$ is homogeneous of degree $-1$ in the variables $a_j$ because $\sigma_i$ has the form $\sum a_j\zeta_j$, and since $\beta(y)=1$, it follows that 
$$(Z_y+\beta(y)) \Pi_\gamma(\tau_0;\sigma)=0.
$$
We have shown that the function $\Pi_\gamma(\tau_0;\sigma)$ is killed by all generators of the tautological system.
\end{proof}

The proof shows that we can relax the condition that $\tau_0$ is an eigenvector of $G$, at the cost of allowing $\tau_0$ to vary in $H^0(L+K_X)$. As a result, the tautological system can be modified to account for this variation. More precisely, we have  

\begin{thm} (Enhanced tautological system) \label{enhanced-tautological}
Let $H-M\ra X$ be a CY bundle and assume {\bf X1-X3}.  Then the period integrals $\Pi_\gamma(\tau;\sigma)$ with $\gamma\in H_{d-s}(Y_\bullet,\Z)$ and $\tau\in H^0(L+K_X)$, satisfy the equations
\begin{eqnarray*}
& &p(\partial_\zeta)\Pi_\gamma(\tau;\sigma)=0\blank(p\in I(\hat X,V))\cr
& &(Z_x+\beta(x))\Pi_\gamma(\tau;\sigma)=\Pi_\gamma(x\cdot\tau;\sigma)\blank(x\in\hat\fg)
\end{eqnarray*}
where $\beta=(0;1,1,...,1)\in{\hat\fg}^*=\fg^*\oplus\C^s$ and $x\cdot\tau=0$ for $x\in\C^s$. 
\end{thm}

\begin{rem}
The last two results can be made slightly more general in two ways. We can replace the full linear system $H^0(X,L_i)$ by a nonzero linear subrepresentation $V_i^*\subset H^0(X,L_i)$ of $G$ in assumption {\bf X1.} We can also allow the linear system to have base locus, in which case, we replace the map $X\ra\PP V_i$ by a rational map $X-\ra\PP V_i$. Then $\hat X$ is given by the closure of the image of the rational map $\phi$ in {\bf X3.} 
\end{rem}

\begin{rem}
The enhanced tautological system in the preceding theorem can be re-interpreted as follows. Fix a basis $\tau_1,..,\tau_m$ of $H^0(L+K_X)$, and form the (column) vector valued sections 
$$\vec\Pi_\gamma(\sigma):=(\Pi_\gamma(\tau_1;\sigma),\cdots,\Pi_\gamma(\tau_m;\sigma))^t.$$
 Then the linear $\hat G$ action on $H^0(L+K_X)$ is realized by a map $\rho_L:\hat G\ra GL_m$, with $(\C^\times)^s\subset\ker\rho_L$. The differential equations corresponding to the symmetry operators then read
$$
(Z_x+\beta(x))\vec\Pi_\gamma(\sigma)=\rho_L(x)\vec\Pi_\gamma(\sigma)\blank(x\in\hat\fg.)
$$
\end{rem}

\bs
\begin{exm}\label{toric-revisited}
CY hypersurfaces in a toric manifold, revisited.
\end{exm}
\cut
Let $X$ be a Fano toric manifold and $T$ the torus acting on $X$ as in Example \ref{toric-example}. We saw that $V^*:=H^0(-K_X)$ has a basis consisting of monomials $z^\alpha=z_1^{a_1}\cdots z_t^{a_t}$ with $\sum a_j[D_j]=[-K_X]$. They correspond to Laurent monomials $y^\mu$, under the correspondence $a_j-1=\bra\mu,\nu_j\ket$ (this is also known as the monomial-divisor map.) Let $z^{\alpha_0},..,z^{\alpha_p}$ be the monomial basis of $V^*$, and $y^{\mu_0},..,y^{\mu_p}$ the corresponding Laurent monomials. Since $-K_X$ is very ample \cite{Oda}, we have a projective embedding afforded by these sections, which by Proposition \ref{function-section}, are represented by the functions $M\ra\CC$, $m\mapsto m^{\alpha_i}$. Thus the embedding reads
$$
X=M/H\into\PP^p,~~[m]\mapsto[m^{\alpha_0},..,m^{\alpha_p}].
$$
The ideal of its image $I(\hat X,V)$ is the binomial ideal $\bra\zeta^{l_+}-\zeta^{l_-}|l\in\cL'\ket$, where $\cL':=\{l\in\ZZ^{p+1}|\sum l_i\bar\mu_i=0\}$ with $\bar\mu_i=(1,\mu_i)$. Here $l_\pm\in\ZZ^p_\geq$ are such that $l=l_+-l_-$. As a special case of Theorem \ref{general-tautological}, with $L=-K_X$, we see that the period integrals 
$$
\Pi_\gamma(1;\sigma)=\int_\gamma\Res \Omega_\sigma
$$ 
of the family $\cY$ of CY hypersurfaces in this case are solutions to the tautological system $\tau(\hat X,V,\hat T,\beta)$ ($\hat T=T\times\C^\times$, $\beta=(0;1)\in{\hat\ft}^*$) which is the GKZ $\cA$-hypergeometric system  with $\cA=\{\bar\mu_0,..,\bar\mu_p\}$ \cite{GKZ1990}\cite{Batyrev}. If we replace $T$ by $\Aut X$, we would recover the extended GKZ system introduced in \cite{HLY1994}.

\section{Examples: homogeneous spaces}

In this section, we apply results we have developed so far, to an arbitrary compact homogeneous space $X$ of a linear algebraic group $G$, to construct period integrals of CY complete intersections in $X$ and to describe their differential equations in details. The special case where $X$ is partial flag variety ($G=SL_n$) was considered in \cite{LSY1}. In this section, we solve the problem in full generality.

First, as an immediate consequence of Theorem \ref{Lichtenstein} and a special case of Theorem \ref{general-tautological}, we have

\begin{thm}\label{homogeneous-one}
Let $G$ be a semisimple group, $X$ be a compact homogeneous space of $G$ viewed as a $G$-manifold, and let $L=-K_X$ and $V=H^0(L)^*$. Then the period integrals 
$$\Pi_\gamma(1;\sigma)=\int_\gamma\R_\sigma(1)
$$
for CY hypersurfaces in $X$ are solutions to the holonomic tautological system $\tau(\hat X,V,\hat G,\beta)$ with $\beta=(0;1)$. The system is generated by the symmetry operators $Z_x+\beta(x)$ ($x\in\hat\fg$), together with the embedding polynomial operators given by the quadrics in Theorem \ref{Lichtenstein}.
\end{thm}

We now give a different description of this system. First, we recall some basic facts about homogeneous spaces. For simplicity, we assume that $G$ is a connected, simply connected, semisimple group. Fix a choice of the simple roots $\Delta$, and let $\{\lambda_\alpha|\alpha\in\Delta\}$ be the dominant integral weights such that $\bra\lambda_\alpha,\beta\ket=\delta_{\alpha\beta}$ for $\alpha,\beta\in\Delta$. Then up to conjugations, parabolic subgroups of $G$ are parameterized by subsets of $\Delta$. Namely, a subset $S\subset\Delta$ corresponds to the parabolic $P_S$ whose Lie algebra is
$$
\fp_S=\fb+\sum_{\beta\in[S]}\fg_{-\beta}
$$
where $\fb$ is the Borel subalgebra of $\fg$ and  $[S]$ is the set of all positive roots (relative to $\Delta$) in the span of $S$. A theorem of Bott's implies that $\Pic(G/P_S)=\sum_{\alpha\in\Delta-S}\ZZ\lambda_{\alpha}$, and that a line bundle $L=\sum n_\alpha\lambda_\alpha$ is ample iff it is very ample iff $n_\alpha>0$ for all $\alpha\in\Delta-S$. 

\begin{prop}\label{KX2}
Put $X=G/P_S$. We have $K_X^{-1}=\sum n_\alpha\lambda_\alpha$ with $n_\alpha\geq2$ for all $\alpha\in\Delta-S$.
\end{prop}
\begin{proof}
Since $\rho:=\half\sum_{\alpha\in\Phi^+}\alpha=\sum_{\alpha\in\Delta}\lambda_\alpha$ and $\bra\lambda_\beta,\alpha\ket=\delta_{\beta\alpha}$, we have $\bra\rho,\beta\ket=1$ for $\beta\in\Delta$. Since $\lambda\equiv K_X^{-1}=G\times_{P_S}\wedge^{top}(\fg/\fp_S)$, we have $\lambda=2\rho-\sum_{\alpha\in[S]\cap\Phi^+}\alpha$. It follows that
$$
\bra\lambda,\beta\ket=2-\sum_{\alpha\in[S]\cap\Phi^+}\bra\alpha,\beta\ket.
$$
For $\alpha$ a nonnegative linear combination of $S$, and $\beta\in\Delta-S$, $\bra\alpha,\beta\ket$ is a sum of nonnegative integers times off diagonal entries of the Cartan matrix for $G$, which are all negative or zero. It follows that the second sum on the right is always nonpositive. With the minus sign in front, it follows that $\bra\lambda,\beta\ket\geq2$.
\end{proof}

Put
$$
\{\beta_1,..,\beta_r\}=\Delta-S
$$
and fix an ample line bundle $L=\sum n_{\beta_i}\lambda_{\beta_i}=\lambda\in \Pic(X)$. We can factor the embedding $\varphi_L:X\into\PP V$, $V=H^0(L)^*$, as follows. Fix an order of $\Delta$ (for $G$ classical simple, we can choose the order to be the natural order given by the flag variety $G/B$, where $B$ is the standard Borel subgroup.) Let $S_i=\Delta-\{\beta_i\}$. Note that the $P_{S_i}$ are exactly the maximal parabolic containing $P_S$. Since $P_S=P_{S_1}\cap\cdots\cap P_{S_r}$, we can define the {\it incidence map}:
$$
\iota: X=G/P_S\into G/P_{S_1}\times\cdots\times G/P_{S_r},~~gP_S\mapsto(gP_{S_1},..,gP_{S_r}).
$$
For each $i$, we have the Pl\"ucker embedding of the generalized Grassmannian 
$$\pi_i:G/P_{S_i}\into\PP W_i,~~gP_S\mapsto[g w_i]$$ 
where $W_i$ is the fundamental representation of highest weight $\lambda_{\beta_i}$ and $w_i$ is a highest weight vector in $W_i$. 
We also have the Veronese maps 
$$\nu_i:\PP W_i\into\PP Sym^{n_{\beta_i}} W_i,~~[v]\mapsto[v\otimes\cdots\otimes v]$$
 and the Segre map 
$$
\psi:\PP Sym^{n_{\beta_1}}W_1\times\cdots\times\PP Sym^{n_{\beta_r}} W_r\into\PP W_L,~~([u_1],...,[u_r])\mapsto[u_1\otimes\cdots\otimes u_r]
$$
where 
$$
W_L:=Sym^{n_{\beta_1}}W_1\otimes\cdots\otimes Sym^{n_{\beta_r}} W_r.
$$
Put 
$$\nu=\nu_1\times\cdots\times\nu_r,~~\pi=\pi_1\times\cdots\times\pi_r.$$
Then we get a $G$-equivariant embedding
$$
\phi_L=\psi\circ\nu\circ\pi\circ\iota:X\into\PP W_L
$$
such that $\phi^*\cO_{\PP W_L}(1)=L$. By the Borel-Weil theorem, $H^0(L)^*=V$ is an irreducible module with highest weight $\lambda$.  Clearly, $\lambda$ is the highest weight in $W_L$ of multiplicity 1, implying that $W_L$ contains a unique copy of $V$.
In particular, there is a unique (up to scalar) $G$-homomorphism
$$
j_L:V=H^0(L)^*\into W_L.
$$
It follows the image $\phi_L(X)$ in $\PP W_L$ lies in the linear subspace defined by $V^\perp\subset W_L^*$, consisting of the linear forms on $W_L$ annihilating $V\subset W_L$. Moreover, $V^\perp$ contains every linear form vanishing on $\phi_L(X)$. 

\begin{defn} (Notations)\label{notations}
Given an ample line bundle $L=\sum n_{\beta_i}\lambda_{\beta_i}=\lambda\in \Pic(X)$, we put
\begin{eqnarray*}
& &W_L:=Sym^{n_{\beta_1}}W_1\otimes\cdots\otimes Sym^{n_{\beta_r}} W_r\cr
& &\phi_L:=\psi\circ\nu\circ\pi\circ\iota:X\into\PP W_L\cr
& &j_L:V=H^0(L)^*\into W_L\cr
& &j_L^*:W_L^*\onto V^*\cr
\end{eqnarray*}
as defined in the preceding paragraph. Let $\hat X$ be the cone over the image of $X\into\PP V$ in $V$. We shall also view $\hat X$ as a $G$-subvariety in $W_L$ via the inclusion map $\hat X\subset V\into W_L$ as well.
\end{defn}

\begin{lem} 
Fix a positive integer $n$ such that $n_i:=n_{\beta_i}\geq n$ for all $i$. Let $X$ be any subvariety of $Y:=\PP W_1\times\cdots\times\PP W_r$. Suppose the vanishing ideal $I(X,Y)$ of $X$ in $Y$ is generated by elements of degrees $(k_1,..,k_r)$ with each $k_i\leq n$.  
Then $I(X,\PP W_L)$ is generated by $I(Y,\PP W_L)$ and the space of linear forms $V^\perp\subset W_L^*$ vanishing on $X$ in $\PP W_L$.
\end{lem}
\begin{proof}
First some notations: let $\cE$ be the set of exponents $v=(v_1,..,v_r)\in\Z^{m_1}_\geq\times\cdots\times\Z^{m_r}_\geq$ such that $|v_i|=\sum_j v_{ij}=n_i$.
Fix a basis $z_{ij}$ ($1\leq j\leq m_i=\dim W_i$) of $W_i^*$. Let $\xi_{i,v_i}$ be the basis of $\Sym^{n_i}W_i^*$ such that $\nu_i^*:\xi_{i,v_i}\mapsto z_i^{v_i}$ (restriction map on sections), and $\zeta_v$ be the basis of $W_L^*$ such that $\psi^*:\zeta_v\mapsto\xi_{1,v_1}\cdots\xi_{r,v_r}$. We denote by $\C[\zeta]_k=\C[W_L]_k$ the subspace of degree $k$ functions on $W_L$, and $\C[z]_{k_1,..,k_r}$ the degree $(k_1,..,k_r)$ functions on $W_1\times\cdots\times W_r$. Put $\phi^*=\nu^*\circ\psi^*:\xi_v\mapsto z^v=z_1^{v_1}\cdots z_r^{v_r}=\prod_{i,j} z_{ij}^{v_{ij}}$.

Obviously that the radical ideal $I(X,\PP W_L)\subset\C[\zeta]$ contains $I(Y,\PP W_L)+\bra V^\perp\ket$. We want to prove the reverse inclusion. So let $p(\zeta)\in I(X,\PP W_L)_m$, $m\geq1$. Since $\C[X]=\C[z]/I(X,Y)\cong\C[\zeta]/I(X,\PP W_L)$, under $\phi^*:\C[\zeta]\ra\C[z]$, we have $\phi^*I(X,\PP W_L)\subset I(X,Y)$, so $\phi^*(p)=p(\zeta_v=z^v)\in I(X,Y)_{n_1m,..,n_rm}$ (each $z^v$ has multi degree $(n_1,..,n_r)$.) By supposition, $I(X,Y)$ is generated by $I(X,Y)_{k_1,..,k_r}$ with all $k_i\leq n\leq n_i$, so
\comment{{\bf 10/20/11 Update.} I notice that for this abstract lemma, the condition $k_i\leq n\leq n_i$ can be weaken to simply $k_i\leq n_i$. Namely, the role of $n$ is completely unnecessary. So, we do not need to know that $-K_X=\sum_i n_i\lambda_i$ with $n_i\geq2$. But in the application below, we ended up needing $n_i\geq2$, because this guarantees $n_i\geq k_i$.}
\begin{eqnarray}
& &I(X,Y)_{n_1m,..,n_rm}=\sum_{k_1,..,k_r\leq n} \C[z]_{n_1m-k_1,..,n_rm-k_r} I(X,Y)_{k_1,..,k_r}\cr
&  &=\C[z]_{n_1(m-1),..,n_r(m-1)}\sum_{k_1,..,k_r\leq n} \C[z]_{n_1-k_1,..,n_r-k_r} I(X,Y)_{k_1,..,k_r}\cr
\end{eqnarray}
Now, restricting the Veronese-Segre isomorphism $\phi^*:\C[\zeta]/I(Y,\PP W_L)\ra\C[z^v:v\in\cE]=Im~\phi^*$ to degree 1 subspace gives an isomorphism (by Borel-Weil for the group $SL(W_1)\times\cdots\times SL(W_r)$, the $z^v$ are linearly independent)
$$
\phi^*:W_L^*=\sum_{v\in\cE}\C\zeta_v\ra\C[z]_{n_1,..,n_r},~~\zeta_v\mapsto z^v.
$$
Since we have $V^\perp=(\phi^*)^{-1} I(X,Y)_{n_1,..,n_r}$ and
$$\sum_{k_1,..,k_r\leq n} \C[z]_{n_1-k_1,..,n_r-k_r} I(X,Y)_{k_1,..,k_r}=I(X,Y)_{n_1,..,n_r}\subset\C[z]_{n_1,..,n_r}$$
it follows that
\begin{eqnarray}
\phi^*(p)\in I(X,Y)_{n_1m,..,n_rm}&=&\C[z]_{n_1(m-1),..,n_r(m-1)}\phi^*(V^\perp)\cr
&=&\phi^*\C[\zeta]_{m-1}\phi^*(V^\perp).
\end{eqnarray}
So $\phi^*(p)\in \phi^*(\bra V^\perp\ket)$. It follows that $p-p'\in ker~\phi^*$ for some $p'\in\bra V^\perp\ket\subset\C[\zeta]$. By Nullstellansatz, $ker~\phi^*=I(Y,\PP W_L)$. This proves that $p\in I(Y,\PP W_L)+\bra V^\perp\ket$.
\end{proof}

We remark that the preceding proof carries over if we use the slightly weaker assumption that $k_i\leq n_i$ for all $i$, hence making the role of the intermediate integer $n$ unnecessary. But for our application below, we do have the integer $n=2$ at our disposal.

\begin{thm}\label{Veronese-quadrics}
Let  $L=\sum n_{\beta_i}\lambda_{\beta_i}=\lambda\in \Pic(X)$ be an ample line bundle such that $n_{\beta_i}\geq2$ for all $i$. Then the ideal $I(\hat X,W_L)$ of $\phi_L(X)$ in $\PP W_L$ is generated by $V^\perp\subset W_L^*$, together with the Veronese binomials 
\begin{equation}\label{e6}
\zeta_u\zeta_v-\zeta_w\zeta_t,\blank u+v=w+t\blank(u,v,w,t\in\cE)
\end{equation}
where the index set $\cE$ and the basis $\zeta_u$ of $W_L^*$ are defined above.
\end{thm}
\begin{proof}
By the theorem of Kostant and Lichtenstein and the Borel-Weil theorem, 
$$(\pi\circ\iota)(G/P_S)\subset\PP W_1\times\cdots\times\PP W_r$$ 
is defined by quadratic forms on the $W_i$, together with multilinear forms on $W_1\times\cdots\times W_r$ (cf. \cite{Kempf1988}.) Since $n_{\beta_i}\geq2$, it follows from the preceding lemma that $I(\hat X,W_L)$ is generated by $V^\perp$, and the vanishing ideal $I(Y,\PP W)$ of $Y\equiv(\psi\circ\nu)(\PP W_1\times\cdots\times\PP W_r)$ in $\PP W_L$.

By the theorem of Kostant and Lichtenstein again, applied to the group $SL(W_1)\times\cdots\times SL(W_r)$ acting on the irreducible module $W_L$, $I(Y,\PP W_L)$ is generated by quadrics. By term-wise elimination, it is straightforward to show any quadric vanishing on $Y$ is a linear combination of the quadrics eqn.\ref{e6}. Conversely, these quadrics obviously vanish on $Y$. This completes the proof.
\end{proof}

We are now ready to give our second description for the tautological system $\tau(\hat X,V,\hat G,\beta)$, which governs the period integrals of the family $\cY\ra B=V^*-D$ of CY hypersurfaces in $X$. The description is based on pulling back the family to $(j_L^*)^{-1}(B)\subset W_L^*$. Notations in Definition \ref{notations} apply. 

\begin{thm} \label{Segre-Veronese}
Let $X=G/P_S$. Put $L=-K_X=\sum_{\alpha\in\Delta-S} n_\alpha\lambda_\alpha$, $V=H^0(L)^*$, and let $W_L,\phi_L,j_L,j_L^*,\hat X$ be as in Definition \ref{notations}.
Then the period integrals $\Pi_\gamma(1;j_L^*\sigma)$ of the family of CY hypersurfaces $\cY$ in $X$, are solutions to the tautological system $\tau(\hat X,W_L,\hat G,\beta)$ with $\beta=(0;1)\in\hat\fg$. This system is generated by
\begin{eqnarray*}
& &Z_x+\beta(x)\blank  (x\in\hat\fg)\cr
& &\partial_{\zeta} \blank(\zeta\in V^\perp\subset W_L^*)\cr
& &\partial_{\zeta_u}\partial_{\zeta_v}-\partial_{\zeta_w}\partial_{\zeta_t}~~\mbox{with}~~u+v=w+t\blank(u,v,w,t\in\cE.)
\end{eqnarray*}
Its solution sheaf is canonically isomorphic to that of $\tau(\hat X,V,\hat G,\beta)$.
\end{thm}
\begin{proof}
The last assertion follows immediately from Lemma \ref{change-of-variables}. Since the canonical isomorphism from the solution sheaf of $\tau(\hat X,V,\hat G,\beta)$ to the solution sheaf of $\tau(\hat X,W_L,\hat G,\beta)$ is given by $f\mapsto f\circ j_L^*$, and since the period integrals $\Pi_\gamma(1;\sigma)$ are solutions to the first system, by Theorem \ref{general-tautological}, it follows that $\Pi_\gamma(1;j_L^*\sigma)$ are solutions to the second system. Finally, Proposition \ref{KX2} implies that for $L=-K_X$, Theorem \ref{Veronese-quadrics} gives the desired generators of $I(\hat X,W_L)$.
\end{proof} 

\begin{rem}
There are at least two ways to explicitly work out the linear forms in $V^\perp\subset I(\hat X,W_L)$. One way is by decomposing the $G$-module $W_L$ into its irreducible component (say by using character theory), and then express elements of $V^\perp\subset W_L^*$, in terms of the basis $\zeta_v$. The second way is by using the theorem of Kostant and Lichtenstein to work out the quadratic forms on each $\PP V_i$, and a result of Kempf's \cite{Kempf1988} to work out the multilinear forms on $W_1\times\cdots\times W_r$, that define the image of the map $\pi\circ\iota$, explicitly. The linear forms they induced on $\PP W_L$ can then be written down easily.
\end{rem}

\begin{rem}
The existence of an ideal like $I(\hat X,W_L)$ that is generated by just linear forms and binomials, has another fairly simple conceptual explanation. It can be thought of as a consequence of the following well-known fact: any projective variety $X$ can be embedded into a sufficient large projective space in such way $X$ can be cut out by linear forms and binomials. Indeed, it follows from the Segre and Veronese embeddings, which are both equivariant embeddings. Therefore, the tautological system $\tau(\hat X,V,\hat G,\beta)$ of any equivariantly projective $G$-variety $X$ can in principle be described by a statement analogous to Theorem \ref{Segre-Veronese}. However, what is special about the case of a homogeneous space is that, the linear forms and binomials in questions can be enumerated quite explicitly as we have done above, thanks to representation theory.
\end{rem}

We now return to the context of \ts\ref{sec-globalPR}, and consider general type complete intersections in a homogeneous space $X$. Since we will be frequently using the tools developed for general CY $H$-bundles, we will denote a parabolic subgroup of $G$ by $H$, rather than $P_S$, in the following discussion.

Fix a linear algebraic group $G$ of dimension $t$, and a parabolic subgroup $H$ of $G$ of dimension $q$. Throughout this section, we put 
$$
M=G,\blank X=G/H,\blank K=G\times H
$$
and view $M$ as a $K$-variety by the usual left and right translations. 

\begin{prop}  \label{CY-structure-GmodP}
Let $\omega_1,..,\omega_t$ be a basis of left $G$-invariant 1-forms on $M=G$, and put
\begin{eqnarray*}
& &\omega_M:=\omega_1\wedge\cdots\wedge\omega_t.\cr
& &\chi_M:H\ra\Aut(\wedge^t\fg^*),~~h\mapsto Ad^*(h^{-1}).
\end{eqnarray*}
Then $(\omega_M,\chi_M)$ is a CY structure on the principal bundle $H-M\ra X$.
\end{prop}
\begin{proof}
Clearly $\omega_M$ is a nowhere vanishing top form on $M$. Note that $\omega_M$ is a basis of the one-dimensional representation $\wedge^t\fg^*$ of $G$, induced by the co-adjoint action. Since $H$ acts on $M$ by right translation, the induced $H$-action on $\wedge^t\fg^*$ is given by $h\mapsto Ad^*(h^{-1})$, as asserted.
\end{proof}

\begin{cor}
If $G$ is semisimple, then $\chi_M=1$ and $K_X\cong L_{\chi_\fh}$.
\end{cor}
\begin{proof}
This follows from the preceding proposition, Theorem \ref{KX-theorem}, and the fact that the $\wedge^t\fg^*$ is a trivial representation of $G$.
\end{proof}

\begin{cor} 
The meromorphic form $\Omega_\sigma$ in Theorem \ref{residue}, up to scalar multiple, is given by 
$$\Omega_\sigma={1\over \sigma_M}\omega_{q+1}\wedge\cdots\wedge\omega_t.$$
\end{cor}
\begin{proof}
Fix a basis $x_1,..,x_t$ of $\fg$ so that $x_1,..,x_q\in\fh\subset\fg$, and let $\omega_1,..,\omega_t$ be the dual basis.
Then the asserted formula for $\Omega_\sigma$ follows from the fact that $\iota_{x_i}\omega_j=\delta_{ij}$. If we change the basis $x_i$, then both $\omega_M=\omega_1\wedge\cdots\wedge\omega_t$ and $x_1\wedge\cdots\wedge x_q$ change by scalar multiples.
\end{proof}

\begin{lem}
Let $G$ be semisimple. Every line bundle on $X=G/H$ is $G$-equivariant and has the form $G\times_H \C_\chi$ for some $H$-character $\chi$. Let $L_1,..,L_s$ be ample and $L+K_X$ is either ample or trivial line bundles on $X$, where $L=\sum_i L_i$. Then {\bf X1-X3} hold.
\end{lem}
\begin{proof}
The first assertion follows from Proposition \ref{Popov}. Since an ample line bundle is very ample on a compact homogeneous space, it follows that the linear systems $H^0(L_i)$ are base point free. By the Borel-Weil-Bott theorem, for $L+K_X$ ample or trivial, $H^0(L+K_X)$ is an irreducible representation of $G$, hence nonzero. So, we conclude that {\bf X1}, hence {\bf X2-X3}, all hold.
\end{proof}

We now let the parabolic subgroup $H$ act on $X=G/H$ by left translations. So, the role of the group $G$ in {\bf X1-X3} is now assumed by $H$ here.

\begin{thm}
Let $X=G/H$ and $L_1,..,L_s,L$ be as in the preceding lemma, and consider the family $\cY$ of complete intersections defined by them as in {\bf X1-X3}. Let $\tau_0$ be a highest weight vector in $H^0(L+K_X)$ with respect to $H$. Then the period integrals $\Pi_\gamma(\tau_0;\sigma)$ of the family $\cY$ are solutions to the tautological system $\tau(\hat X,V,\hat H,\beta)$ with $\beta:=(\lambda_0;1,1,..,1)$, where $\lambda_0\in\fh^*$ is the eigenvalue corresponding to $\tau_0$. Moreover, the system is holonomic.
\end{thm}
\begin{proof}
By Proposition \ref{CY-structure-GmodP}, the principal bundle $H-G\ra X$ admits a CY structure. By the preceding lemma, {\bf X1-X3} hold. By Theorem \ref{general-tautological}, our first assertion follows. Let $B$ be the Borel subgroup of $G$ contained in $H$. Then it is well-known that $G/B$ has only a finite number of $B$ orbits under the left translation action. Since $G/B\onto G/H$, $X=G/H$ also has only a finite number of $B$ orbits. Hence the same is true under the $H$ action. It follows that the cone $\hat X$ has only a finite number of $\hat H$ orbit. Now the last assertion follows from Theorem \ref{holonomic}.
\end{proof}

\section{Concluding remarks}

To further understand period integrals using the framework developed in this paper, one immediate problem is to compute their power series expansions. As we have indicated earlier, this can be done, at least in some special cases, by using the global Poincar\'e residue formula. For example, in the case of CY complete intersections in a toric manifold, explicit power series solution can be obtained by manipulating the global meromorphic form (see for example \cite{HLY1996}) inside the dense orbit of the torus, and then integrate it against certain real torus cycle. One can try to do the same for homogeneous spaces. One way to do so is by first picking an affine chart on $X$ that is a copy of $\C^d$, and then attempt to integrate the global meromorphic form against a similar real torus cycle in $\C^d$. In the toric case, explicit power series solutions near the so-called maximal unipotent point in the parameter space, are also known in great generality \cite{HLY1994}. They could be used as a guide for obtaining similar formula in the case of homogeneous spaces. 

In the algebraic setting, the entire framework we have developed for the global Poincar\'e residue map in this paper carry over to any characteristic zero varieties. In positive characteristics, most of the results carry over as well. However, some new phenomena occur. For example, the  the period integrals of a complete intersection family in {\bf X1-X3} are still solutions to a tautological system $\tau(\hat X,V,\hat G,\beta)$. However, the proof of Theorem \ref{general-tautological} clearly shows that the ideal of polynomial operators that annihilates the period integrals is much larger than that the embedding ideal $I(\hat X,V)$. In fact,  if the field of definition has charateristic $p$, then any polynomial $Q\in\C[V]$ of degrees $(k_1,..,k_s)$ such that $k_i\geq p$ for some $i$, gives a new differential operator $Q(\partial_\zeta)$ that annihilates the period integrals. This may be an indication that the enlarged tautological system is much more restrictive on the period integrals than in the characteristic zero case. Both problems mentioned here will be investigated in a future paper.


\address B.H. Lian, Department of Mathematics, Brandeis University, Waltham MA 02454. lian@brandeis.edu.

\address 
S-T. Yau, Department of Mathematics, Harvard University, Cambridge MA 02138. 
yau@math.harvard.edu.

\end{document}


\maketitle



The purpose of this erratum is to correct `Theorem 3.12 (Uniqueness of CY structure)' and its proof in the paper {\it Period integrals of CY and general type complete intersections}, {\sl Invent. Math. Vol. 191, 1 (2013) 35-89.} The correct statement is


\bs
\noi{\bf Theorem 3.12} {\it Let $H-M\ra X$ be a $G$-equivariant principal $H$-bundle over $X$. Then we have the 1:1 correspondence}
\[
\{\text{$G$-equivariant CY structures on $M$} \} 
\longleftrightarrow
 \lambda_{M,G}^{-1}([K_X]_G),~~(\C\omega,\chi)\lra\chi\chi_\gh  
 \]
{\it where} $\lambda_{M,G}:\Hom(H,\C^\times)\ra\Pic_G(X)$, $\chi\mapsto [M\times_H\C_\chi]_G$.

\bs
Here $\chi_\gh$ is the $H$-character for the one-dimensional $H$-module $\wedge^{\dim\gh}\gh$ induced by $\Ad:H\ra\Aut\gh$.
This shows that $G$-equivariant CY structures on $M$ need not be unique in general, but they are classified precisely by the coset
$$
\chi\chi_\gh+\ker\lambda_{M,G}\subset\Hom(H,\C^\times)
$$
in the $H$-character group, if $(\C\omega,\chi)$ is a given $G$-invariant CY structure on $M$ (Theorem 3.3 \cite{LY}). The rest of the paper \cite{LY} is not affected, except for the following errata: on the twenty fourth line on page 39: and prove that when exists, a CY structure is unique iff the natural map $\Hom(H,\C^\times)\ra\Pic(X)$ is injective (Theorem 3.12); on the twentieth line on page 46: holomorphic $G$-invariant top form; right side of eqn. (2.3) on p45: $([m],\mu_\alpha(m))$; Proposition 3.13: character $\chi\chi_\gh$; line two of p53: 
$$mh^{-1}\mapsto[mh^{-1},1]=[m,(\chi\chi_\gh)(h)^{-1}]=[m,1](\chi\chi_\gh)(h)^{-1}.
$$
Line 16 on f p53:
$$
\Omega_\sigma:={1\over\sigma_M}\Omega\in H^0(\Omega^d_M|(M-\{\sigma_M=0\})).
$$

\bs\noi{\it Proof of Theorem 3.12.}
Let $\lambda_{M,G}(\chi')=[K_X]_G$. By Theorem 3.3\cite{LY}, there is a $G$-invariant CY structure of the form $(\C\omega,\chi'\chi_\gh^{-1})$ on $M$ which is mapped to $\chi'$ under the given correspondence. This shows that this correspondence is onto. If $(\C\omega_1,\chi_1)$ and $(\C\omega_2,\chi_2)$ are $G$-equivariant CY structures mapped to the same character, then $\chi_1=\chi_2$. By Lemma 3.9\cite{LY}, $\C\omega_1=\C\omega_2$. This shows that the correspondence is bijective. $\Box$ 



\bs
Here is an example \cite{CL} where non-uniqueness can occur.

\bs
\noi{\bf Example.} Put $G=SL_{n+1}$, $H=(\C^\times)^2$, $X=\P^n$, $M=\cO(-1)^\times\oplus\cO(-1)^\times$, the Whitney sum of two copies of the complement of the zero section in $\cO(-1)$. Then $\lambda_{M,G}:\Hom(H,\C^\times)\simeq\Z^2\ra\Pic_G(X)=\Pic(X)=\Z$, $\chi_{k,l}\mapsto k+l$, where $\chi_{k,l}(h_1,h_2)=h_1^kh_2^l$. Since $[K_X]=-n-1$, a complete list of CY structures on $M$ is given by
$$
(\C\omega_l,\chi_{n+1-l,l}),~~l\in\Z
$$
where $\omega_l$ is the unique, up to scalar, nonvanishing top form on $M$ which transforms by the $H$-character $\chi_{n+1-l,l}$. Explicit formulas for the forms $\omega_l$ are given in \cite{CL}. 
